\documentclass[11pt]{article}
\usepackage{AVDocLarge}

\usepackage{graphics}
\usepackage{amssymb}
\usepackage{color}
\begin{document}


\title{
Individual-based models\\
under various time-scales
 }
\author{Aurelien Velleret
\footnote{Aix-Marseille Universite, CNRS, Centrale Marseille, I2M, UMR 7373
	13453  Marseille, France, 
	email~: aurelien.velleret@univ-amu.fr}
}

\maketitle

\section*{Abstract}
This article is a presentation
of specific recent results describing 
scaling limits of individual-based models.
Thanks to them, 
we wish to relate the time-scales
typical of demographic dynamics
and natural selection
to the parameters of the individual-based models.
Although these results are by no means exhaustive,
both on the mathematical and the biological level,
they complement each other.
Indeed, they provide a viewpoint for many classical time-scales.
Namely, they encompass 
the time-scale typical of the life-expectancy 
of a single individual,
the longer one wherein 
a population can be characterized
through its demographic dynamics,
and at least four interconnected ones 
wherein selection occurs.
The limiting behavior is generally deterministic.
Yet, since there are  selective effects 
on randomness
in the history of lineages,
probability theory
is shown to be a key 
factor in understanding the results.
Besides,
randomness can be maintained 
in the limiting dynamics,
for instance to model
rare mutations fixing in the population.

 Les r\'{e}sultats r\'{e}cents pr\'{e}sent\'{e}s
 dans cet article 
 d\'{e}crivent des limites d'\'{e}chelles 
 de  mod\`{e}les individus-centr\'{e}s.
Gr\^{a}ce \`{a} eux, 
nous allons mettre en lumi\`{e}re les diff\'{e}rentes \'{e}chelles de temps caract\'{e}ristiques 
des variations d\'{e}mographiques 
et de la s\'{e}lection naturelle.
L'objectif n'est pas d'\^{e}tre exhaustif,
tant au niveau math\'{e}matique que biologique.
N\'{e}anmoins, 
ces r\'{e}sultats sont tr\`{e}s compl\'{e}\-mentaires
les uns des autres.
Ils fournissent un aper\c{c}u des principales \'{e}chelles de temps d'int\'{e}r\^{e}t.
Ils englobent notamment
l'\'{e}chelle de temps de la vie d'un individu,
celle plus longue o\`{u} la population
 peut \^{e}tre d\'{e}crite comme une entit\'{e}
avec ses caract\'{e}ristiques propres
qui dirigent les dynamiques
d\'{e}mographiques,
et au moins quatre autres \'{e}chelles 
imbriquées
sur lesquelles la s\'{e}lection joue un r\^{o}le.
La dynamique limite est g\'{e}n\'{e}ralement d\'{e}terministe.
Pour autant,
puisque l'action de la s\'{e}lection 
se base sur l'al\'{e}a pr\'{e}sent dans l'histoire des individus,
les probabilit\'{e}s apparaissent comme un \'{e}l\'{e}ment cl\'{e} pour comprendre ces résultats.
Par ailleurs, 
la stochasticit\'e peut aussi \^etre conserv\'ee \`a la limite,
par exemple pour mod\'eliser 
de rares \'ev\'enements de fixation de mutations.

\section*{Introduction}

The origin of this paper
is the session "Stochastic Processes and Biology" during the conference Journ\'ees MAS in Dijon at the end of August 2018.
Four talks were presented, 
with a variety of methods and probabilistic approaches.
Their common feature, that motivates this presentation,
is the fact that they start from similar Individual-Based Models
 (IBM)
 and justify an asymptotic behavior as 
 both the population size and the observation time are large.
 The time-scales 
wherein these observations can be obtained are nonetheless different,
 yet interconnected.
The structure of the paper
aims  both
at presenting the principal objects 
of the four talks of the session
and at reflecting this increasing sequence of time-scales involved
in adaptation of populations.
\\
 
For more than fifteen years,
there has been a significant activity 
around a new probabilistic framework 
for mathematical modeling of ecology, 
population genetics 
and trait evolution.
Traits are any features of the individuals
whose transmission the modeler is interested in,
and we refer to the beginning of Chapter 1 
for some classical examples.
They may vary
during the life-time of the individuals,
at births
or only at rare mutation events.
The population is described 
by the distribution of its traits,
from which one may 
for instance deduce its growth rate.
Any effect of natural selection
should depend on these traits
and affect their distribution.

The IBM
are a priori the best formulation
for validation 
of a macroscopic model of population dynamics.
These processes detail at individual level 
births,  deaths and interactions 
in the population.
The description may also specify
migration patterns, 
aging 
or competition for resources 
between individuals.
This setting is clearly the closest to  actual simulations 
designed by computational biologists
to validate theoretical models,
or simply some predictions.
In this view, 
such setting is also the one 
with the least number of simplifying assumptions.
It is also certainly the one where calculations 
of probabilities of events
are the most complicated.
One of the first aims of such probabilistic modeling 
has actually been 
to connect these individual-based models and simulations
to much simpler systems of Ordinary Differential Equations (ODE)
or Partial Differential Equation (PDE) models.

To describe natural selection
on heritable traits,
systems of ODE
usually provide the most classical models 
to simply express their effects.
They have been refined 
by introducing PDE
to deal with a continuum of traits 
in populations.
Thanks to these models,
one can cover most 
of actual mathematical models
of ecology, population genetics 
and character evolution.
Yet, this deterministic approach 
does not take into account 
the variability observed 
when reproducing similar experiments.
A part of this variability is due to external perturbations.
Stochastic Differential Equations (SDE)
have been introduced 
to model these effects.
They usually assume 
that the parameters appearing in the systems of ODE / PDE
are themselves subject to stochasticity.
The other main component of this noise is internal 
and inherent to the discreteness of the system :
selection between traits emerges
 through competition
between a finite number of individuals.
Taking this component into account 
is the purpose of these IBM.
\\

To justify simpler models from the IBM,
the usual assumption
is that populations are large 
enough so that the law of large numbers 
makes random fluctuations negligible.
Given the finite population sizes 
ecologists are interested in,
an elementary description may well be accurately designed
while a finer structuring of the population is much more demanding.


It is rather natural for mathematicians
to try to describe simplified behavior through an asymptotic 
on large time-scales 
(i.e. for large number of generations).
It is 
	at the core of Darwin's theory of natural selection
to assume a separation of time-scale 
	between apparently incidental randomness 
		among contemporary individuals
	and the evolution of traits, 
		and even of species.
The goal is surely 
to obtain exact results of convergence 
	to a non-trivial process 
		for the new dynamics in a new time-scale.
It is however very challenging, 
and perhaps an approximation cannot be justified
in this way,
while still being helpful.
There is also no reason to restrict ourselves 
purely to deterministic limiting processes :
a first order of fluctuations
can be inferred from central limit theorems.
For instance, it may  take the form 
of a solution to a specific SDE.
Such a next order
 notably helps to specify the conditions of validity 
of the deterministic approximation.
	
An archetypal example
is the stochastic approximation of the population size
in a continuous setting (i.e. with Brownian-type fluctuations).
Beyond such a result is the idea 
	that individual events 
	(of births or deaths)
	have quite a negligible impact 
	on the population as a whole,
yet randomness still occurs 
through varying accumulation of such effects.
Thus, to obtain the convergence
towards a stochastic continuous process,
one needs (besides the population size going to infinity)
to amplify 
the frequency of births and deaths also.
Notably, this  means  that selective effects appear 
in a much longer time-scale than the life-expectancy 
of individuals 
(cf. in particular chapter 3 in \cite{BM15}).
The scale of population size 
provides the essential relation between the parameters
at individual time-scale 
(on this time-scale, the birth rate is around 1)
and at population level 
(cf. section \ref{sec:aline}).
In case there is more than one parameter to adjust,
analysis may be even less intuitive 
and may even lead to different behaviors associated 
with different scalings.
\\

Each of the four talks given 
during the above-mentioned session constitutes 
the foundation of each of the four sections
to follow.
The results that we will present
are taken from the referenced articles,
and the reader is spared some details 
on the underlying assumptions.
These choices 
and the discussions that complement these results
express my personal view on the subject,
and do not engage the authors of the articles.
My aim while reporting 
these results is to highlight the variety 
of time-scales involved in population models.
These population models are simpler to analyze
and quicker to simulate 
than the IBM for large populations.
We will thus specify 
how large the population sizes shall be 
in order to have valid estimations
of  individual-based models
through population models.

In Section \ref{sec:aline},  
we shall present the work of Hoffmann and Marguet \cite{Aline}.
IBM are introduced 
in terms of traits that structure the population and evolve
in this section during the life-time of individuals.
The estimation of this trait dynamics 
is the purpose of their paper,
in which  specific estimators
are proposed and their accuracy evaluated.
From this evaluation,
we will obtain some insight 
into the time-scale required for different populations 
(with different trait dynamics)
to be distinguished.
The following sections are devoted to evolution 
and the associated models rely on traits
that are much more accurately transmitted to descendants.
Although this previously described variability is usually neglected in evolutionary studies
(only some average effect is considered),
we shall see in Section \ref{sec:aure} 
that it may play a significant role
when two components of selection 
induce competing effects.
Both the strength 
and the time-scale of selection
between two sub-populations 
may in fact involve the amount of a priori non-selective variability.
This is one of the main conclusions 
that I wished to discuss based on my work in \cite{V19}.
Assuming however that this variability has only 
a minor effect on the 'local' reproductive value 
of individuals
and rescaling time properly,
 Champagnat and Henry describe in \cite{Henry}
the process of evolution
 as an almost deterministic behavior.
 This work and the connections
 to the two previous sections
 will be our concern in Section \ref{sec:henry}.
 Finally, 
 natural selection may
 rely on the emergence of very rare mutations,
as in Tran's talk.
 The associated time-scale 
 depends on the occurrence rate 
and the mean effect of successful mutations
 (that invade the population).
Among several possibilities,
I have chosen to present 
the models from  
\cite{bfmt18}
and \cite{BBC17}
in Section \ref{sec:tran}.

\section{Recognition of population traits from individual dynamics}
\label{sec:aline}

This section is devoted 
to the results of Hoffmann and Marguet.
For simplicity, 
the model neglect any interaction between individuals
other than transmission of specified traits from parents to offspring,
what we call a model of branching population.
Assuming that changing traits
induce a structure 
on this branching population,
they propose a statistical estimation 
of this dynamics of traits \cite{Aline}.
This is the right place 
to discuss this notion of individual-based models 
and how to relate individual histories 
to the law that governs these behaviors.
It gives some insight into the number of generations involved 
before one can actually distinguish distinct population models.
On a shorter time-scale,
natural selection may favor some trajectories
by mere chance.
Yet those lucky behaviors
should be observed with close frequencies 
between distinct populations.
This selection shall thus hardly lead 
to any heritable effect.
For instance, as presented  in \cite{Aline},
any statistical estimation of traits
requires a sufficiently large population history.
Likewise, 
a long-term history is required 
for natural selection to favor 
one subpopulation 
over another (with different inheritance).
Moreover,
an initially very favorable combination 
of trait and environment
might not be so beneficial in the long term :
most descendants 
might be exposed later on
to a much less favorable environment
given their traits.
This "sample size" shall thus be larger 
when there are strong correlations
between successive generations.
In cases 
where such correlations are sufficiently weak,
Hoffmann and Marguet (cf. \cite{Aline})
study how to take them into account in statistical analysis.
They indicate the accuracy of some estimators 
in terms of the size of the genealogical tree
on which observations are indexed.
The focus in \cite{Aline}
is on the random process governing birth events,
for which they propose the first statistical analysis 
in a structured population.
\\

At an individual-based level,
each individual is characterized by some value $x\in \mathcal{X}$,
that we call a trait.
Typically $\mathcal{X} \subset \mathbb{R}^k$ or $\mathcal{X} \subset \mathbb{R}^k\times F$,
where we may consider $k$ characteristics
 such as size, 
 spatial position,
 age,
the amount of certain proteins,
of certain ressources,
of parasites and so on,
while $F$ denotes a finite (or numerable) set of classes,
such as sex,
eye color etc.
Before the death of the individual,
this trait $x$ may be constant,
but it can also evolve 
according to some stochastic differential equation (SDE).
Interaction with other individuals 
can then appear 
in parameters of this equation,
possibly depending on the trait $x'$ of the other individuals around.
Death of individuals 
happens at a rate $d$,
that may depend first on the trait $x$ of the individual,
but also on the whole population
and eventually its effect on the environment
(as is the case in Section 3).
We recall that we could include 
the age of the individual 
in the information carried by $x$, 
so that there is almost no restriction.
Death at a fixed age is also not difficult to include in the model.
Same kind of dependencies 
may be considered for births,
and apply to the number of offspring 
and their states at birth.
Independent exponential variables
or Poisson Point Processes are thus usually
a very convenient way to encode all of these events.
In a nutshell, 
one can a priori 
represent with stochastic IBM models
any computational model that one could design
for validation of biological predictions,
as soon as one can follow the individuals one by one.
I refer notably 
to Subsection \ref{sec:iBm} for a more concrete description
of such measure-valued process.
In the following, 
we will focus the analysis
on continuous trait-spaces,
for which one can exploit the regularity 
of the estimated functions.

\subsection{Discussion on the assumptions in \cite{Aline}}

\subsubsection{A process on an incomplete tree}

The main specificity of this statistical approach
is to deal with the dependency 
given by the genealogy of individuals.
For instance, 
looking at the size as a trait,
the larger the mother-cell is,
the larger its daughter-cells shall be.
With the experiments 
where one follows each cell individually,
one can also easily obtain 
the entire associated lineages.
Yet, 
many of these cells are no longer observable 
because of the design
of experimental processes.
Although the authors do not mention it,
death of some cells is likely to be included
as well.
The statistical analysis of the traits
thus relies on data 
that is indexed by some tree 
that is "incomplete" as compared 
to the complete binary tree.
A more detailed definition 
on the assumption on the tree
is given below.
In order to relate the accuracy of the estimators
to the number of birth events, 
the authors consider
both the cases of a bounded population size
at each generation
and the case of an population size
expanding at a  given growth rate~:
\begin{Def}
Consider the complete binary tree, 
on which our tree will be indexed.
Define the generation $n$ of this complete tree as $G_n$.
Note that $Card(G_n) = 2^n$.

A regular incomplete tree is a family of subsets $U_n$
of the complete tree up to generation $n$ such that~:

$\quad$ (i) the parent of any individual of $U_n$ is also in $U_n$
 
 and  (ii)  for some $0 \le \rho \le 1$~:\\
$ 0 < \liminf_{n\rightarrow \infty} 2^{-\rho n} Card (U_n \cap G_n) 
\le \limsup_{n\rightarrow \infty} 2^{-\rho n} Card (U_n \cap G_n) <\infty$.
\end{Def}

In the case $\rho \in (0,1]$, 
the data characterize 
a growing population of cells.
For instance, it corresponds 
 to an experiment 
where one lets these cells duplicate 
freely in a rich medium.
In the case $\rho = 0$, 
one rather expects some population size at equilibrium,
although it is not required here.
It also corresponds to the experiments 
where only one cell is followed after each division
(which is the case in recent experiments 
in a microfluidic environment : the mother-machine).

In the following results, 
this indexation tree
is considered as given.
It means that we consider 
the trait dynamics
conditionally on the realization of the tree.
Yet,
 the interest of the experimenter
is a priori rather 
on intrinsic parameters without conditioning.
Thus,
there is an implicit assumption 
that the shape of the indexation tree 
is independent 
of the trait dynamics.

With a single lineage experiment, 
this is not really an issue.
But if one has to consider natural death events
while estimating variable birth rates,
assuming the independence
 seems much more questionable.
It appears a bit surprising
to have a birth rate 
depending on the trait
while the probability of giving birth 
in one life-span
is independent of it.
A constant death rate would be 
a more natural assumption
than a prescribed tree.
Yet, it would be much more difficult to analyze
and would possibly not reflect 
the design of the experiment  very well.

\subsubsection{Ergodicity of the trait}

The other main assumption in \cite{Aline}
is the very rapid convergence to the ergodic measure $\nu$
in the transmission of the trait at birth
 from the parent to its offspring.
Assuming that $u\in U_{n}\cap G_n$ 
describes the parental index 
of an individual with index
$v \in U_{n+1}\cap G_{n+1}$,
the trait at birth $X_v$ of the latter 
is given by $Q(X_u, dx_v)$
where $X_u$ is the trait at birth of the parent.
The convergence of the estimators 
is stated uniformly
as long as  there exist $C > 0$ 
and a positive weight-function $V$ 
such that for any $ m \ge 1$, 
 and $\varphi$ satisfying 
 $\|\varphi\|_V := \sup_{x} |\varphi(x)|/[1+V(x) ]<\infty$~:
 \begin{align}
\|Q^m(\varphi) - \nu(\varphi)\|_V \le  C\, 2^{-m}\, \|\varphi - \nu(\varphi)\|_V.
\label{Q12}
 \end{align}
The set of such  $Q$ 
is denoted $\mathbb{Q}_{1/2}$.
 Note in this view that the authors 
 choose exploit exclusively 
traits at birth 
 for their statistical estimations.
 
\subsubsection{Definition of the stochastic process, regularity and confinement}

The model in \cite{Aline}
is concerned with a trait 
supposedly governed during the life-time 
of each individual
 by a stochastic flow
of the form :
$
dX_t = b(X_t)\, dt + \sigma(X_t) dW_t,\quad
$
where $W$ refers to a standard Brownian Motion,
and $(X_t)$ evolves on $\mathcal{X} \subset \mathbb{R}^k$.
Birth occurs at time $s$ at rate $B(X_s)$.
At reproduction event, 
the authors assume that the value of the trait 
is distributed between two offspring
with the following mechanism :
given an independent r.v. $\theta$,
the trait at birth of the two offspring is given
respectively by $\theta y$ and $(1-\theta) y$.
$\theta$ is  drawn according to $\kappa(y) dy$,
 for some probability density
 function $\kappa(y)$ on $[0, 1]$.

This model can for instance 
represent the size of the cell
or the propagation of a parasite
(that we take here as an archetype). 
The main issue  is here 
to estimate the effect of the amount of parasites in a cell
on its reproduction. 
Note that the individual level is 
the one of a given cell, and not of a given parasite, 
as it would be in a completely discretized model.
Yet, the number of parasites is assumed to be so large
that a continuous description by a random process 
is a well-justified simplification.
The estimation is focused on the effect of these parasites
on the birth rates $(B(x) : x\in \mathcal{X})$
 of the hosts
and on the law of distribution 
$(\kappa(y) : y\in [0,1])$ between offspring.
\\

Assumptions are 
specifically designed in \cite{Aline}
for such a model.
In a broader perspective,
we rather focus on the core principles 
for which they are introduced
and refer to \cite{Aline} 
for the precise statements.
By Assumption 2 on the drift and diffusion coefficient,
the authors ensure 
that the trait $\mathcal{X}$ stays somewhat confined around $0$,
with a uniformly elliptic and non-singular diffusion.
By Assumption 3 on the birth rate,
they ensure some regularity in $x$, 
a boundedness condition 
on the potential explosion of births
and they prevent the vanishing of these events.
Finally, by Assumption 4 
on the splitting of $x$ at birth between the two newborns,
they ensure a lower-bounded density on the fragmentation parameter,
and prevent too asymmetric partitioning.
The authors mention
that this assumption could probably be relaxed.

\subsection{Main results for the estimations
 of the Generation kernel
and the birth rate}

The first step of the analysis 
consists in the estimation of the kernel $Q$
and its stationary distribution $\nu$.
This step should be easily generalized
in a broader perspective 
of kernels $Q$. 

For any function $\psi$ : $\mathcal{X}\times \mathcal{X}\rightarrow \mathbb{R}$,
and $y\in \mathcal{X}$,
consider $\psi_{\star}(y) := \sup_{x\in \mathcal{X}}|\psi(x, y)|$,
$\psi^{\star}(x) := \sup_{y\in \mathcal{X}}|\psi(x, y)|$.
Denoting by $\wedge$ the minimum,
we also define
for any positive measure $\rho$ 
on $\mathcal{X}$~:
\begin{align*}
&|\psi|_\rho 
:= \int_{\mathcal{X}^2} |\psi(x, y)|\, \rho(dx)\, dy
+ \left ( \int_{\mathcal{X}^2} |\psi(x, y)|\, dx\, dy
\wedge \int_{\mathcal{X}} |\psi_{\star}(y)|\, dy\right )
\\&\hcm{3}
\mathcal{M}_{\mathrm{U}_{n}}(\psi)
:= \dfrac{1}{Card(\mathrm{U}_{n}^\star)}\, \sum_{u\in \mathrm{U}_{n}^\star} \psi(X_{u-}, X_u),
\end{align*}
where  $\mathrm{U}_{n}^\star$ is 
$\mathrm{U}_{n}$ deprived from the root,
and $X_{u-}, X_u$ denotes respectively the trait at birth of the parent of $u$ 
and the one of individual $u$ itself.
Recall that the transition from $X_{u-}$ 
to $X_u$ is given by $Q(X_{u-}, dx)$.

\begin{prop}
{\it Let Assumptions 2, 3 and 4 be satisfied. 
Let} $\mu$ {\it be a probability measure on} $\mathcal{X}$ {\it such}
{\it that} $\mu(V^{2}) < \infty$. 
{\it Let} $\psi$ : $\mathcal{X}\times \mathcal{X}\rightarrow \mathbb{R}$
 {\it a bounded function such that}
  $\psi_{\star}$ {\it is compactly supported. 
  If}
$\mathrm{U}_{n}$ {\it is a regular incomplete tree, the following estimate holds true}:
\begin{align}
\mathbb{E}_{\mu} [(\mathcal{M}_{\mathrm{U}_{n}}(\psi) - \nu(Q \psi))^{2}]\
\lesssim Card(\mathrm{U}_{n})^{-1} (|\psi^{2}|_{\mu+\nu}+|\psi^{\star}\psi|_{\mu}
+(1+\mu(V^{2}))|\psi_{\star}|_{1} |\psi|_{\nu}),
\label{CVM}
\end{align}
{\it where the symbol} $\lesssim$ {\it means up to an explicitly computable constant that depends only on supp} $(\psi_{\star})$ {\it as long as} $\mathrm{Q}\in \mathbb{Q}_{1/2}.$
More generally, it would also depend on $Q$.
\end{prop}

For the following theorem,
we assume that the operator $Q$ 
has a density $(q(x, y))_{(x, y) \in \mathcal{X}^2}$ w.r.t. Lebesgue measure on $\mathcal{X}$,
with some H\"{o}lder regularity.
I refer to \cite{Aline}
for the exact definition of the sets 
$\mathbb{Q}_{1/2}^{\alpha,\beta}(R)$,
where the value $R>0$ defines a bound 
on the H\"{o}lder regularity 
of respectively order $\alpha$ in $x$ and $\beta$ in $y$. 
The "$1/2$" refers to the property (\ref{Q12}).
Given such a regularity,
the authors propose to adjust the order of the estimation kernel 
and explain how to choose the associated window sizes 
given $|\mathrm{U}_n|$ (for some large $n\ge 1$).
The estimation of $q$ also depends on a threshold $\varpi_{n}$
that is to adjust.
I refer to \cite{Aline} 
for the exact definitions of the estimators 
$\hat{\nu}_{n}(y)$ 
and $\hat{q}_{n}(x, y)$
of respectively $\nu(x)$ and $q(x, y)$
(with the knowledge of $\mathrm{U}_n$).

\begin{theo}
\label{nuQ}
{\it Let Assumptions 2, 3 and 4 be satisfied. 
	Assume that the initial distribution} $\mu$ {\it is}
{\it absolutely continuous w.r.t. the Lebesgue measure with a locally bounded density and}
{\it satisfies} $\mu(V^{2}) <\infty.$
{\it Let} $\alpha, \beta > 0$.
{\it Then, for any $\varrho$-regular incomplete tree} $\mathrm{U}_{n}$ 
{\it  and any} $R>0,$
\begin{align*}
&\hcm{4}
\underset{Q\in \mathbb{Q}_{1/2}^{\alpha,\beta}(R)}{\sup}\
(\mathbb{E}_{\mu}[(\hat{\nu}_{n}(y)-\nu(y))^{2}])
^{1/2}
\lesssim Card(\mathrm{U}_{n})^{-\beta/(2\beta+1)}.
\\&
\text{and}\hcm{2}
\underset{Q\in \mathbb{Q}_{1/2}^{\alpha,\beta}(R)}{\sup}\
 (\mathbb{E}_{\mu}
 [(\hat{q}_{n}(x, y)   - q(x, y))^{2}])^{1/2}
 \lesssim \varpi_{n}^{-1}
 Card(\mathrm{U}_{n})^{-s(\alpha,\beta)/(2s(\alpha,\beta)+1)}
\end{align*}
{\it hold true, where} 
$s(\alpha,\beta)^{-1} := (\alpha\wedge\beta)^{-1}+\beta^{-1}$
 {\it is the effective anisotropic smoothness associated}
{\it with}~$(\alpha,\beta)$.
\end{theo}
Note that the rate of convergence depends on the regularity of $q$.
The higher it is, 
the more useful information we can gather 
from observed transitions from vicinities of $x$ 
to vicinities of $y$.
The accuracy  of the estimators is thus better supported.
\\

The estimation of the division rate is much more complicated.
Additional assumptions are required.
The study is restricted to diffusion of the trait
 on a compact space, 
 so that (\ref{CVM}) can be strengthened with a right-hand side
 only depending on the $L^\infty$-norm of $\psi$.
 Again, we restrict the study to $\mathrm{Q}\in \mathbb{Q}_{1/2}.$
A parametric approach is considered, 
where $B$ is
encoded by parameters $\vartheta \in \Theta$.
The associated Fisher matrix $\Psi(\vartheta)$
 is assumed to be non-singular
 to avoid issues of identifiabilty.
 They also assume some upper-bounds 
 on the derivative of $B$
 along $\vartheta$ up to the third order.
Moreover, 
the authors prevent any degeneracy and singularities 
of the birth rate.
Finally, this birth rate $B$ is assumed 
to be a globally monotone function
of  $\vartheta$
(uniformly in $\mathcal{X}$).

\begin{theo}
\label{theta}
 {\it Let the above-mentioned assumptions be satisfied. For every} $\vartheta$ {\it in the interior of}
$\Theta$, {\it if} $\mathrm{U}_{n}$ {\it is a} $\varrho$-{\it regular incomplete tree :}
$$
\sqrt{Card(\mathrm{U}_{n})} \times 
(\hat{\vartheta}_{n}-\vartheta)
\rightarrow \mathrm{N}(0, \Psi(\vartheta)^{-1})
$$
{\it in distribution as} $n \rightarrow \infty$, {\it where} $\mathrm{N}(0, \Psi(\vartheta)^{-1})$ {\it denote the} $d$-{\it dimensional Gaussian distribution}
{\it with mean} $0$ {\it and covariance the inverse of} $\Psi(\vartheta)$.
\end{theo}

Note that we find in this case the classical rate of convergence 
with order the square root of the number of observations,
as for classical Markov Chains.
It hints at the fact that every observed transition from parents to offspring
contributes to the estimation of $\vartheta$.

The monotonicity condition of $B$ along $\vartheta$
seems very restrictive,
especially in the case 
of multidimensional parameter.
Relaxing this condition 
would be a natural direction to look at for improvement.
\\

\textbf{Concluding remarks :}

As mentioned after the assumption 
on the genealogical tree,
considering this tree as prescribed
is much more convenient.
The fact that the convergence results 
of Theorems \ref{nuQ} and \ref{theta}
 depend on $U_n$ only through its cardinal 
 indicates that these results 
 shall be quite robust 
 concerning the specific realization of the tree.
 
However, it does not exclude 
that the estimators are slightly biased
by the fact that certain transmission patterns 
 might lead to a greater probability of survival 
 (of the lineages).
For instance, 
it may be beneficial to divide early, 
in order not to let the parasites divide for long 
before the next division.
Traits favoring early division 
are thus more likely to be transmitted.
This bias depends on how strongly these effects 
may have shaped the genealogical tree.
Given that the design of the experimental process
has a strong effect on the shape of the tree,
including this dependency 
would though be very difficult 
and likely to introduce even more bias.
Besides, 
the estimation of the biased birth rate 
may be the main interest. 
It only means that one infers the dynamics 
of a "successful" lineage,
with luck playing a role in this success,
as well as environmental conditions. 
\\

This issue of estimating individual parameters 
by looking at the population as a whole 
leads to similar issues of independence.
Notably,
the population around the focal individual
can be considered as a component of its environment.
One may wish that the parameters of its dynamics
depend on the interaction it has with this extended environment.
The models become much simpler if one can neglect
 the detailed interactions between the individuals,
 and replace them by average effects.
 In the limit of a large population, 
 for a system at equilibrium,
 such assumption may be justified by the Law of Large Number.
 Indeed, the effects of the interactions are averaged,
 as long as they are not too local, 
 and thus almost globally constant.
 In results of propagation of chaos,
 such property is generalized 
to cases where the population is not at equilibrium.
 I refer to Sznitman's lecture at Saint-Flour 
 \cite{S91}
 for an overview on this topic.
 In such a limiting case, 
 known as McKean-Vlasov equations,
 the law of the process itself acts 
 on the individual dynamics,
 together with fluctuations specific to this individual.
 At equilibrium, 
 this law is all the more stable
because the population size is large.
 Thus, we can deal with the effects of these interactions 
 through some hidden parameters 
 (like an effective death or birth rate).
And as for the analysis of this section,
estimating them by simply 
looking at death and reproduction times
would produce a bias if they depend on heritable factors.
But this bias is presumably small 
if heritability is weak.
 \\

To study evolution of population traits,
that is for instance $\vartheta$ if mutations could alter this value,
it is very classical to assume that the population size is large.
Since one usually assumes a separation of time-scales
between demographic dynamics 
and evolutionary processes,
there is clearly time 
for individual fluctuations 
to be averaged.
Heritable traits 
can leave a sufficiently clear mark
for natural selection to be effective.
When the environment itself 
is affected by the traits of individuals,
as in the next model,
it might however not be so clear what an average effect would be.

\section{Selection with two levels }
\label{sec:aure}

The foundation of this section
is my work on the ability 
of selective effects 
acting at a group level 
to compensate for those acting inside each group
at an individual level \cite{V19}.
For simplicity,
we focus here 
on competition between two types 
in a population of fixed and large size.
In real populations and IBM models, 
there are fluctuations in the proportion 
due to the inexact compensation
between births and deaths events 
occurring with the same large rate.
Biologists usually refer to it 
as genetic drift,
and often neglect its effect
in the case of large populations.
One shall see however that 
these fluctuations may not be neutral at all 
in a model 
where two selective effects are considered :
the first one favors some individuals inside their groups
(individual level selection)
while the other favors some groups 
depending on the individuals they gather
(group level selection).
Reducing random fluctuations
inside each group 
strongly hinders
response to selection at group level.
For clarity,
we shall assume that these selective effects
 are conflicting.
For instance,
one may ask
if and how 
the inefficient or cheater individuals
can be regulated through natural selection 
at this group level.
\\

The individual-based model is taken from \cite{LM15},
where a formalism for group selection is introduced.
All groups have the same size $n \in \mathbb{N}$.
There are two types of individuals : $C$ and $D$.
Type $D$ individuals
have a better reproduction
at the individual level 
($D$ for defectors) 
while type $C$ individuals 
are positively selected 
at the group level
($C$ for cooperators). 
Replication and selection occur concurrently 
at individual and group levels 
according to a "nested Moran process",
as introduced in \cite{Dur08}
and recalled next. 
 Type $C$ individuals replicate at rate $w_I$
  and type $D$ individuals 
  at rate $w_I\,(1 + s), s \ge 0$. 
When an individual gives birth, 
another individual in the same
group is selected uniformly at random 
to die and be replaced,
so that the population size remains constant. 
To reflect antagonism at the higher level of selection, 
groups replicate at a rate 
that depends on  the number of type $C$ individuals they contain.
We take this rate 
to be $w_G\times \,[1 + r(k/n)]$, 
where $k/n$ is the fraction 
of  type $C$ individuals in the group
and $r(x), x\in [0,1]$ 
is the selection coefficient
 at group level.
The number of groups 
is maintained at $m$ by selecting a group uniformly at random 
to die whenever a group replicates. 
The two offspring of groups are assumed to be identical to their parent.
\\

Two limits of large population 
are justified from the individual-based model 
in \cite{LM15}.
At least four different effects may 
contribute to the limiting dynamics.
The selective effect at individual level
is due to difference 
in growth rate between $I$ and $G$ individuals;
the selective effect at group level
to difference 
in the growth rate of groups
(depending on the proportion of $G$ individuals).
Consider the number of replacement events 
in a Moran model
of a population 
made up of a constant number of identical individuals.
In a limit of large populations 
as stated in the central limit theorem,
this number
shall increase with a linear rate
with random fluctuations
that may be well 
approximated by a Brownian Motion.
Similar fluctuation effects generate
variations in the proportion of trait carriers 
in a population.
One can include such effects,
that we call "random fluctuations",
in the limiting large population process,
both at individual 
and group levels.
%

Both limits in \cite{LM15}
 include a term for each selective force, 
but only  the second includes
 random fluctuations,
both within the groups and between groups.
In \cite{V19}, 
another limit is considered,
where random fluctuations  
are kept  only inside the groups,
together with selective forces.
Although such limit can be mathematically justified,
it is undoubtedly more realistic to keep also
random fluctuations at group level. 
Yet, the analysis is notably much more technical
with much less clear interpretation.
In this more complete model,
one type go almost surely 
extinct in finite time.
Such event of "ultimate fixation"
has a positive probability for each type.
For large enough group populations,
the limit obtained by neglecting these fluctuations 
provides an interesting view 
on the main features of the dynamics.
\\

Let $X^i_t$
 be the number of type $C$ individuals in group $i$ at time $t$. 
 Then~:
\begin{equation*}
\textstyle 
\mu^{m;n}_t 
:= \frac{1}{m} \sum_{i\le m} \delta_{X^i_t/n}
\end{equation*}
is the empirical measure at time $t$ 
of the proportion of type $C$ by group,
with $m$ the number of groups  
and $n$ the number of individuals per group. 
Here, $\delta_x$ is the Dirac at $x$. 
The $X^i_t$ are divided by $n$ so that $\mu^{m;n}_t$  is a
probability measure on $E_n := [0;1/n; ... ; 1].$
For fixed $T > 0$, $(\mu^{m;n}_t)_{t\le T} \in D([0; T ]; \mathcal{M}_1(E_n))$, 
the set of c\`{a}dl\`{a}g processes on $[0; T ]$
 taking values in $\mathcal{M}_1(E_n)$
 (the set of probability measures on $E_n$).
 With the particle
process described above, $\mu^{m;n}_t$ has generator
\begin{equation*}
(\mathcal{L}^{m;n} \psi )(v) 
= \textstyle \sum_{i,j}
(w_I\, R_I^{i,j} + w_G\, R_G^{i,j})(v)
\times \left(  
\psi\big [v+ 1/m\, (\delta_{j/n} - \delta_{i/n})\big] - \psi[v]\right)
\end{equation*}
where $\psi \in C_b(\mathcal{M}_1([0; 1]))$ is a bounded continuous function, 
and $v \in \mathcal{M}_1(E_n) \subset \mathcal{M}_1([0; 1]).$ 

The transition rates $(w_I\, R_I^{i,j} + w_G\,R_G^{i,j})$ are given by
\begin{align}
&&R_I^{i,j} (v) 
&:= \left\{
\begin{aligned}
&m\, v(i/n) \, i\, (1 - i/n)\, (1+s) 
&\text{if } j = i - 1; i < n,
\\& m\, v(i/n) \, i\, (1 - i/n) 
&\text{if } j = i + 1; i > 0,
\\& 0 &\text{otherwise}\hcm{1.2}
\end{aligned}
\right.&
\end{align}
\begin{equation}
\text{ and }\qquad	 
R_G^{i,j}  (v) 
:=  m\, v(i/n) \,v(j/n)\, (1 + r[j/n] ).
\label{RI}
\end{equation}
$R_I^{i,j}$
and $R_G^{i,j}$ 
are the rates of 
respectively 
individual-
and group-level reproductive events.

\begin{theo}
\label{MLim}
Suppose that  $w_I/n \rightarrow  \omega_I$, 
$n\, s \rightarrow \sigma$
as $n, m\rightarrow \infty$,
while $w_G$ and $\{r(x)\}_{x\in [0,1]}$
are kept constant.
Suppose the particles in the process $\mu^{m;n}_t$ 
are initially independently and 
identically distributed according to the measure $\mu^{m;n}_0$,
 where $\mu^{m;n}_0 \rightarrow \mu_0$ as $m, n\rightarrow \infty$.
 Then, 
$\mu^{m;n}_t$ converges weakly to 
$\mu_t \in D(\mathbb{R}_+; \mathcal{M}_1([0;1]))$,
where $\mu_t$ is the unique solution 
with initial condition $\mu_0$ of 
the differential equation~:
\begin{align}
&\qquad \partial_t\, \langle\mu_t\, \big| \, f\rangle
 = \langle \mu_t\, \big| \, \mathcal{L}_{W\!F}f \rangle
 +  \langle \mu_t\, \big| \, r\,f\rangle  
 - \langle \mu_t\, \big| \, f\rangle 
 \times \langle \mu_t\, \big| \, r\rangle,
 \label{eqCar}
  \\& \text{ where }	
  \mathcal{L}_{W\!F} f(x) = - s\,x\,(1 - x)\, \partial_{x}f(x)
  +(\sigma^2/_2)\,.\, x\,(1 - x)\, \partial^2_{xx} f(x).
\label{LWF}
\end{align}
\end{theo}

\subsection{Definition as a conditional law}
We describe the solution of such equation 
through a Feynman-Kac penalization 
of the stochastic process with generator $ \mathcal{L}_{W\!F}$.
Let us define $X$ as
 the solution of the SDE :
 \begin{align}
 &dX_t := -s\, X_t\, (1-X_t) \, dt 
 + \sigma\, \sqrt{ X_t\, (1-X_t)} \, dB_t,
 \qquad X_0\sim \mu_0
 \label{WFdiff}
 \end{align}
 The existence and uniqueness of such process can be found e.g. in chapter 5.3.1 of \cite{D10}.
It is linked there with an individual-based model, 
namely the neutral 2-allele Wright-Fischer Markov Chains.
 We also consider the 
 following Feynman-Kac penalization~:
\begin{equation*} 
 Z_t := \textstyle{\exp\int_0^{t} r(X_s)\,ds},
  \end{equation*}
 Note that $r$ is bounded so that for any $t>0$,
$\mathbb{E}(Z_t)\in (0,\infty)$.

\begin{prop}
With the above definitions,
we characterize $\mu_t$ by the fact that
for any $f\in \mathcal{C}^2_b$ :
\begin{equation*}
 \langle \mu_t\, \big| \, f\rangle 
 := \mathbb{E}\left[ f(X_t)\, Z_t \right] 
 /\, \mathbb{E}\left[ Z_t \right], 
 \end{equation*}
\end{prop}

I refer to \cite{V19} for the proof of the proposition, 
which is easily derived from the Ito formula.
\\

We next give another view on this law $\mu_t$,
that will simplify our notations and clarify our point.
Since subtracting a constant to $r$ does not change the value of 
 $\langle \mu_t\, \big| \, r\,f\rangle  
  - \langle \mu_t\, \big| \, f\rangle \; \langle \mu_t\, \big| \, r\rangle$,
  we assume in the following that $r\le 0$, 
  so that we can consider it as a death rate.
$Z_t$ can then be interpreted as the probability
that the process has survived until time $t$,
while confronted to a death rate of $r$,
conditionally on $(X_u)_{u\ge 0}$. 
More formally, 
with $T_\partial$ 
an independent exponential r.v.
with mean $1$,
we can define the extinction time as~:
\begin{align*}
&\tau_{\partial} := \inf\left \lbrace t\ge 0\, ; \, 
-\ln(Z_t) \ge T_\partial \right \rbrace,
\quad \text{so that }
\mathbb{P}(t < \tau_{\partial}
\, \big| \, X_u, u>0)
=  \mathbb{P}\left( -\ln(Z_t) < T_\partial
\, \big| \, X_u, u\le t \right) 
= Z_t
\\&\hcm{3}
\langle \mu_t\, \big| \, f\rangle  
= \mathbb{E}\left[ f(X_t)\, ; \, t < \tau_{\partial} \right] 
/\, \mathbb{P}\left[t < \tau_{\partial}\right]
= \mathbb{E}(f(X_t) \, \big| \, t < \tau_{\partial} ).
\end{align*}
Any equilibrium of $\mu_t$ is thus a quasi-stationary distribution (QSD)
of $X$ under the death rate $r$.
Applying recent results on QSD
makes it possible 
to obtain the limiting behavior 
and the speed of convergence of $\mu_t$.
\\

Similarly,
the hitting times of $0$ and $1$ 
are denoted $\tau_0$ and $\tau_1$.
Since they are absorbing, 
it is natural to be interested in the law of the marginal
with these extended extinction times~:
\begin{align*}
\tau_{0, \partial} := \tau_{\partial} \wedge \tau_0, 
\quad \tau_{1, \partial} := \tau_{\partial} \wedge \tau_1, 
\quad \tau_{0,1} := \tau_0 \wedge \tau_1,
\quad \tau_{0,1, \partial} := \tau_{\partial} \wedge \tau_0 \wedge \tau_1.
\end{align*}
\begin{prop}
	With the above definitions,
	we can also characterize $\mu_t$ 
	by the fact that
	for any $f\in \mathcal{C}^2_b$ :
$$\mu_t = x^0_t\, \delta_0 + x^1_t\, \delta_1
+x^\xi_t\, \xi_t$$
	\begin{align*}
&\text{ where }	 x^\xi_t := 
\dfrac{ \mathbb{E}\left[ Z_t\, ; \, t < \tau_{0,1} \right]}{\mathbb{E}\left[ Z_t \right]}, \qquad
\langle \xi_t\, \big| \, f\rangle 
:= \dfrac{\mathbb{E}\left[ f(X_t)\, Z_t\, ; \, t < \tau_{0,1} \right] }
{\mathbb{E}\left[ Z_t\, ; \, t < \tau_{0,1} \right]}
= \mathbb{E}\left[ f(X_t)\, \vert \, t < \tau_{0,1,\partial} \right],
\\& \hcm{1}
x^0_t := \dfrac{\mathbb{E}\left[ Z_{\tau_0}\,\exp[-r_0(t-\tau_0)]\, ; \, \tau_0 < t \right]}{\mathbb{E}\left[ Z_t \right]}, \qquad 
x^1_t := \dfrac{\mathbb{E}\left[ Z_{\tau_1}\,\exp[-r_1(t-\tau_1)]\, ; \, \tau_1 < t \right]}{\mathbb{E}\left[ Z_t \right]}.%
	\end{align*}
\end{prop}


The proof is elementary and left to the reader.
\\

Remarks :
$\bullet$
It is not difficult to generalize the model
to include frequency dependent effects of
selection at individual level.
We would then replace $s>0$ by some smooth function $(s(x))_{x\in [0,1]}$.
Generalizations of $\sigma$ as a function 
are also not a mathematical issue.
Yet a priori, it does not seem  biologically justified.

$\bullet$
This Feynman-Kac penalization
is analogous to the many-to-one formula
described by \cite{BDMT11}
in the case of 
Markov processes associated to branching Galton-Watson Trees.
 Behind this description
  is the idea that there is a bias
towards larger offspring in the genealogy 
of a typical individual in the population,
favoring the individuals with a larger birth rate.
Contrary to \cite{BDMT11}, 
where no interaction occurs between lineages,
one cannot expect 
to represent the law of $\mu_t$
directly through another Markov Process,
i.e. without penalization
(note in (\ref{eqCar}) the quadratic term in $\mu_t$).

\subsection{QSDs and exponential convergence}
Note first that in any case, 
$\delta_0$ and $\delta_1$ are QSDs for the extinction time $\tau_{\partial}$, 
i.e. stable distributions 
for the dynamics given by (\ref{eqCar}).
If the initial condition $\mu$
supported on $\{0, 1\}$,
the dynamics is immediately deduced 
from the death rates in $0$ and $1$.

For other initial conditions $\mu$ (not supported on $\{0, 1\}$),
we define the following semi-groups associated to our different extinctions ~: 
\begin{align*}
&
\mu A_t (dx) 
:= \mathbb{P}_\mu (X_t \in dx \, \big| \, t <\tau_{\partial}),
\hcm{.5}
\mu A^{01}_t (dx) 
:= \mathbb{P}_\mu (X_t \in dx \, \big| \, t <\tau_{0,1, \partial}),
\hcm{0.5}
\mu A^{1}_t (dx) 
:= \mathbb{P}_\mu (X_t \in dx \, \big| \, t <\tau_{1, \partial})
\end{align*} 

\begin{prop}
\label{p:QSD01}
There exists a unique QSD $\alpha\in \mathcal{M}_1[(0,1)]$ 
and a survival capacity of $\eta$
associated to the extinction time $\tau_{0,1, \partial}$.
With the associated extinction rate $\rho_\alpha$, it means~:
\begin{align*}
\frlq{t>0}
&\mathbb{P}_\alpha (X_t \in dx\, ; \, t< \tau_{0,1,\partial}) 
= \exp[-\rho_\alpha\, t]\, \alpha(dx),
\quad
\frl{x}
\eta(x) = \exp[\rho_\alpha\, t]\, \mathbb{E}_x (\eta(X_t) \, ; \, t< \tau_{\partial}) 
\end{align*}
Moreover, we have the following exponential convergences at rate $\zeta>0$  :
\begin{align}
&\hcm{2}
\Ex{C>0}
\frlq{\mu \in \mathcal{M}_1[(0,1)]}
\left\|  \mu A^{01}_t
 - \alpha \right\|_{TV}
\le C \, \exp[-\zeta\, t].
\label{CVal}
\\&\hcm{2}
\Ex{C'>0} \frlq{\mu \in \mathcal{M}_1[(0,1)]}
|\exp[\rho_\alpha\, t]\, \mathbb{P}_\mu(t<\tau_{0,1, \partial})
- \langle \mu\, \big| \, \eta\rangle | 
\le C' \, \exp[-\zeta\, t]
\label{CVeta}
\\&\text{ a fortiori }
\eta(x) := \lim_{t\rightarrow \infty} \exp[\rho_\alpha\, t]\,
\mathbb{P}_x(t<\tau_{0,1, \partial})
\quad\text{ and }
\|\eta_\bullet\| 
:= \underset{\{x\in (0,1),\, t>0\}}{\sup}
\exp[\rho_\alpha\, t]\, \mathbb{P}_x(t<\tau_{0,1, \partial}) < \infty
\label{etaB}
\end{align}
\end{prop}
Let $\rho_0 = - r_0$ 
($\rho_1 = - r_1$)
the extinction rate of $\delta_0$
(resp. $\delta_1$).
We show in the following that the long-time behavior
of the process with only the local extinction rate
depends mainly on $\rho_\alpha,$ $\rho_0$ and $\rho_1$.

In the following convergences, 
we will often have uniform bounds for probability measures 
belonging for some $n \ge 1$ and $\xi >0$ to~:
\begin{align*}
&\mathcal{M}_{n,\, \xi}^{0} 
:= \left \lbrace \mu\in \mathcal{M}_1([0,1])\, \big| \,
\mu[1/n, 1] \ge \xi\right \rbrace\, , \,
\textstyle{\bigcup_{n, \xi}} \mathcal{M}_{n,\, \xi}^{0}  
= \mathcal{M}_1([0,1]) \setminus \{\delta_0\}.
\\\text{or in}\quad
&\mathcal{M}_{n,\, \xi}^{0,1} 
:= \left \lbrace \mu\in \mathcal{M}_1([0,1])\, \big| \,
\mu[1/n,\, 1-1/n] \ge \xi\right \rbrace\, , \,
\quad(n\ge 3, \xi >0)
\\&\hcm{1}
\textstyle{\bigcup_{n, \xi}} \mathcal{M}_{n,\, \xi}^{0,1}  
= \mathcal{M}_1([0,1]) \setminus 
\{u\, \delta_0 + (1-u)\, \delta_1\, \big| \, u \in [0,1]\}.
\end{align*}

\subsubsection{Despite the fixation events,
polymorphic groups maintain themselves in the population}

\begin{prop}
\label{p:rar10}
Assume that $\rho_\alpha < \rho_0\wedge \rho_1 := \rho$.
Then, there is only one stable QSD $\alpha_{0,1}$, 
with convergence rate $\rho -\rho_\alpha$,
i.e. :
\begin{align*}
\frl{n\ge 1}\frl{\xi >0}
\Exq{C_{n, \xi}>0}
\frlq{\mu \in \mathcal{M}_{n,\, \xi}^{0,1}}
\left\|  \mu A_t - \alpha_{0,1} \right\|_{TV}
\le C_{n, \xi} \, \exp[-(\rho_\alpha -\rho)\, t],
\end{align*}
where $\alpha_{0,1}$ has extinction rate $\rho_\alpha$
and is given as 
$\alpha_{0,1} = y_0\, \delta_0 +  y_1\,\delta_1  + y_\alpha\, \alpha$
 with~:
 \begin{align*}
&\dfrac{y_0}{y_\alpha} 
 = \dfrac{\rho_\alpha\times \mathbb{P}_\alpha(\tau_0 = \tau_{0,1, \partial}) }{
  (\rho_0 -\rho_\alpha)}\, , \,
  \quad
 \dfrac{y_1}{y_\alpha} 
  = \dfrac{\rho_\alpha\times \mathbb{P}_\alpha(\tau_1 = \tau_{0,1, \partial}) }{
   (\rho_1 -\rho_\alpha)}\, , \,
 \end{align*}
 and of course $y_0 + y_1 + y_\alpha = 1$.

\end{prop}

If $\rho_1 < \rho_0$, for any initial condition 
$\mu_0 =u\, \delta_0 + (1-u)\, \delta_1$ with $u\in(0,1)$,
$\mu_t$ converges at rate $\rho_0 - \rho_1$ to $\delta_1$.

If $\rho_1 = \rho_0$, then any such distribution is a QSD
 with the extinction rate $\rho_0$.
\\

Pure groups 
are continuously generated
from polymorphic groups
without any reversed transition.
Yet, in this case,
these polymorphic groups 
are sufficiently selected upon
through their better survival 
to persist in the population.
Pure groups are like remnants 
of these polymorphic groups :
their proportion reaches a steady state 
where their faster decay
compensate for the fixation rate
from polymorphic groups.
The stabilization of polymorphic profile 
induces the convergence of both this fixation rate
and of the maintenance rate of polymorphic groups
 to $\rho_\alpha$.

In any case,
polymorphism is maintained 
by any sufficiently large group selection favoring it, since~:
\begin{prop}
\label{prop:rI}
Given any $\sigma >0$, $s\ge 0$, 
and a bounded continuous function $r^0$ with its maximum 
only in the interior of $(0,1)$,
there exists a critical value $R_\vee >0$ such that
for any $R>R_\vee$ and considering the system with $r = R\, r^0$,
we indeed have $\rho_\alpha < \rho_0\wedge \rho_1$.
\end{prop}

Conversely, when group selection is too small, 
polymorphic groups cannot be maintained :
\begin{prop}
\label{prop:r0}
Conversely,
given any $\sigma >0$, $s\ge 0$, 
and a bounded measurable function $r^0$,
there exists a critical value $R_\wedge >0$ such that
for any $R<R_\wedge$ and considering the system with $r = R\, r^0$,
it holds $\rho_0\wedge \rho_1< \rho_\alpha$.
\end{prop}

Too strong neutral fluctuations
also make the fixation of the groups
hardly avoidable, so that~:
\begin{prop}
\label{prop:SI}
Given any $s>0$ and any bounded function $r$, $\lim_{\sigma \rightarrow \infty} \rho_\alpha(\sigma) = + \infty$.
\end{prop}

\paragraph{\textsl{Conditions for polymorphism to be maintained :\\} }
For this maintenance rate $\rho_\alpha$
to be higher than the decay of pure groups,
it is clearly necessary that $r$ is maximal in $(0, 1)$.
Of course, it is not sufficient,
because first of random fluctuations ($\sigma$)
and second because of selection effects inside each group.
As stated in Proposition 2.6, 
too large $\sigma$
 would induce a too large rate of fixation
 (to 0 or 1).
 Even strong effects of selection through $r$ 
 would then be unable 
 to make their reproduction large enough 
 to compensate for this loss.
 
 But on the other hand, 
 when there are internal selective effects 
 pushing towards type $I$ individuals,
having only small random fluctuations 
 limits the effectiveness 
 of selection at group level.
 Indeed, 
 all groups with similar initial condition
 evolve too closely for such effects of selection 
 to really distinguish between them :
 they are essentially driven by the flow of the ODE :
 $\quad 
 \partial_t x_t := -s\, x_t\, (1-x_t).\quad $
 Even if perturbations are amplified by this selection
towards the opposite direction,
strong deviations would be much too costly.

\subsubsection{Fixation on either side is the most stable case 
and the type $C$ is favored by group selection}

\begin{prop}
\label{p:r1r0ra}
Assume that $\rho_1 < \rho_0 < \rho_\alpha$.
Then, $\delta_1$ is the only stable QSD, 
with convergence rate $\rho_0 -\rho_1$~:
\begin{align*}
\frl{n\ge 1}\frl{\xi >0}
\Exq{C_{n, \xi}>0}
\frlq{\mu \in \mathcal{M}_{n,\, \xi}^{0}}
\left\|  \mu A_t - \delta_1 \right\|_{TV}
\le C_{n, \xi} \, \exp[-(\rho_0 -\rho_1)\, t].
\end{align*}
\end{prop}

We also have an additional level of convergence :
\begin{prop}
\label{p2:r1r0ra}
Assume that $\rho_0 < \rho_\alpha$. Then, there exists $C>0$ s.t.~:
\begin{align*}
\frlq{\mu \in \mathcal{M}_1([0,1])\setminus\{\delta_1\}}
\left\| \mu A^{1}_t
 - \delta_0  \right\|_{TV}
\le C \, \exp[-(\rho_\alpha -\rho_0)\, t].
\end{align*}
\end{prop}

Both results are asymptotic 
and might not reflect exactly the dynamics on a short time-scale.
Yet, their justification gives us some insight 
into what can happen.
\\

\paragraph{\textsl{When polymorphism gets quickly negligible :\\} }
If $\sigma$ is large,
except for initial conditions very close to 0,
a non-negligible proportion of pure type $C$ groups 
quickly emerges and dominates the distribution.
For initial conditions very close to 0,
the emergence time of these pure $C$ groups
mainly depends on the proportion of groups 
able to quickly escape such vicinity of 0.
Soon, the growth 
in the proportion of pure $C$ groups
is essentially due to the difference in growth rate 
between these groups and the rest of the population,
meaning that, at this time,
 the fixation 
of polymorphic groups
plays a negligible role.
Finally, one observes 
the competition between the two types of pure groups,
with the initially rare $G$ groups
outnumbering the first dominant $I$.
\\

\paragraph{\textsl{When trajectories are drifted with little fluctuations :\\} }
As explained in the previous subsection,
the flow of the equation 
$\partial_t x_t := -s\, x_t\, (1-x_t)$
dominates the dynamics as long as $\mu_t$ stays localized. 
If $\sigma$ is rather small,
the profile $\mu_t$ is essentially 
given by the integration of the growth rate 
along the trajectories of the flow.
Notably, 
consider for simplicity
the case where the initial condition
is supported on $[0, 1-2 \epsilon]$.
The pure flow brings 
the group from a proportion $1-\epsilon$
to $\epsilon$ in a deterministic time $t_\epsilon$.
If $\sigma$ is small enough,
then with a probability close to 1
the process is upper-bounded 
by this deterministic flow 
starting from $1-2 \epsilon$.
It implies that
$\mu_{t_\epsilon}$ is mainly concentrated 
in $[0, 2\epsilon]$.
If the proportion of pure $G$ groups is still negligible
at this time, 
we shall observe a decay 
in the proportion of polymorphic groups
larger than an exponential 
rate of $\rho_\alpha - \rho_0$.
It is a priori unclear that the rate $\rho_\alpha - \rho_0$
is actually observed,
since the QSD $\alpha$ might be localized around $1$
and very difficult to reach 
for initial conditions with many less 
cooperative groups.
One may expect some stabilization to occur,
where $\mu_t$ restricted to $(0,1)$ gets close to some $\tilde{\alpha}$.
This distribution $\tilde{\alpha}$ is presumably supported mainly close to 0,
with a much larger extinction rate $\rho_{\tilde{\alpha}}$ 
than $\rho_{\alpha}$ (and $\rho_0$).
One shall have the right intuition 
by considering $\tilde{\alpha}$ 
(resp. $\rho_{\tilde{\alpha}}$) 
instead of $\alpha$ (resp. $\rho_\alpha$)
in the reasoning of Propositions \ref{p:r1r0ra} and \ref{p2:r1r0ra}.
This view is supported by first hints of simulations (not detailed in the article).

As stated in Proposition \ref{p:r1r0ra},
pure $G$ groups shall prevail even in that case.
A rate $\rho_1- \rho_0$ of decay 
of remaining groups 
is very likely to be seen,
but possibly after a domination by type I groups
for a large period of time.
The duration of this domination
 depends actually much on the initial distribution,
and particularly on its tail near 1.
Indeed, transitions to 1 is especially costly
in term of its probability of occurrence.
\\

\paragraph{\textsl{Compensate the flow of invasion by $I$ individuals :\\} }
Section 3 provides 
an evaluation
of the strength of this group selection
needed to compensate the flow
in the limit of vanishing $\sigma$.
In this large deviation regime, 
the process seems to evolve
for most of the time
 according to a modification
of the initial flow.
Yet, it is not as simple :
 for instance,
one might observe the abrupt emergence 
of a type 
which was so far negligible,
whose growth rate is much higher 
than the previously dominant type.

This is exactly what shall presumably happen
in this model in the case $r_1< r_0$.
Undoubtedly,
the initial proportion of $G$ groups
or the neutral fixation 
for initially almost pure $G$ groups
concerns a very tiny proportion of ancestors.
Yet,
the ancestors initially drifted
towards larger proportion of $I$ individuals
soon loose any $G$ individual (with very few exceptions).
So
there is a point in time at which 
we can no longer neglect 
the more prolific descendants 
the former exceptional ancestor groups
will have.
This shall certainly correspond to the time 
at which $G$ individuals eventually dominate.
Are the ones that finally dominate 
necessarily pure groups ?
Due to the potential long-term persistence of the process
in the vicinity of $1$ 
when the random fluctuations are very small,
the behavior of $\rho_\alpha$ as $\sigma$ tends to 0
is quite unclear. 
So it might happen that some almost pure $G$ groups actually dominate.
\\

\paragraph{\textsl{Selection upon the initial condition :\\} }
Still,
the selection between 
different initial conditions 
may be effective 
if those are sufficiently apart.
It may postpone for some significant time 
the trend towards 0.
Yet, 
for the polymorphism to persist for long,
a very specific form
of the law of the initial condition
is required.
This has been specified in \cite{LM15} 
in the limit where the random fluctuations are neglected.
Such long-term persistence is effectively possible 
because the flow is vanishing in the vicinity of 1.
The authors consider 
only functions $r$ that are linear 
in the proportion of type $G$ individuals.
But it should generalize to a much more general cases,
provided $r$ is larger near 1 than near 0.

Note however that in such regime,
the groups dominating at a given time 
are essentially not the ancestors 
of those dominating at a much larger time.
The descendants of the former have  
been transported 
by the flow towards 0,
where this sub-population has decayed faster.
On the other hand, 
the ancestors of the latter
must have stayed for long 
in a very tiny and specific region very close to 1,
and are thus very few (while the former dominate).
Indeed, going backward through the ancestry lines
means essentially following the flow backwards.
This backward flow goes quickly towards 1,
which it approaches at exponential speed.
With initial conditions irregular in the vicinity of 1,
we thus might observe a surprising sequence 
of vanishing and reemerging polymorphism
between periods of domination by $I$ individuals.
\\

\paragraph{\textsl{Is this observable in the individual-based model ?\\} }
Any transition involving a reasonable amount of groups
is expected to be indeed observed.
This can be estimated through the number of ancestors 
from the initial population
upon which such transitions shall rely.
The fluctuations shall vanish with the number of them.
Provided the neutral fluctuations
in the births and deaths of groups
are small,
the approximation should be qualitatively valid
with approximately 20 ancestries.

It means also that too exceptional transitions 
are very unlikely to be observed.
For instance,
the escape from a too close vicinity of 0 
happens  with a too small probability.
The most likely is to observe the complete fixation.
For the fixation of pure $C$ groups to be observed,
the most probable is then to have one group
escaping the vicinity of 0, reaching the other boundary
and generating a sufficiently large family there
 for the extinction to become negligible.
Only after such exceptional realization 
becomes the fixation of pure $C$ groups likely to occur.
Given Theorem \ref{MLim},
larger population sizes 
makes the event more likely to occur.
Yet, 
in order for one group 
to behave in a way so different 
from the typical one,
it might be required that
the population size is 
at a largely unrealistic level.
Similarly, 
when $\sigma$ is so small that transitions towards 1 
become negligible, 
and for an initial condition with a light tail (towards 1),
the fixation of pure $C$ groups happens after a very exceptional behavior.

This mathematical complexity is presumably not so relevant 
in terms of the biology.
As soon as the initial condition has a sufficiently light tail
in this vicinity of 1,
one mainly observe a massive proportion of the groups fixing 
as pure $I$ types
and becoming dominant for a very long time.

\subsubsection{Polymorphic groups are more stable than 
type I groups but less than pure type $C$ groups
}

This case is also treated in \cite{V19}. 
The result is a combination of the ones 
in the two previous subsections.
The dynamics for the domination by pure $G$ groups
relies on similar principles as in Proposition \ref{p:r1r0ra}.
The main difference is
that the intermediate convergence
is stated
towards a polymorphic QSD rather than pure $I$ groups.
This polymorphic QSD $\alpha_1$
is described as in Proposition \ref{p:rar10}
when one subtracts pure $C$ groups before the renormalization.
The asymptotic rate of convergence 
towards the Dirac at pure $G$ groups
is deduced from this intermediate convergence result :
$\rho_1-\rho_\alpha$ (smaller than $\rho_1 - \rho_0$).

Again,
for the approximated IBM,
 one may be faced to the same limitations regarding the origin 
 of the first pure $G$ groups as in Proposition \ref{p:r1r0ra}.
In practice, 
the convergence towards the polymorphic QSD
might also not actually reflect the main dynamics of convergence.
It might happen that the "emergence" of pure type 1 groups
can actually be almost concomitant 
to the emergence of $\alpha^1$.
Looking at simulations 
for small values of $\sigma$,
an alternative metastable regime around 0
may dominate for a significant time 
the marginal law restricted to $(0, 1)$.
Comparing with the case of initial condition close to 1,
this distribution $\tilde{\alpha}$ is very different from the actual QSD $\alpha$,
with very separate supports.
Even if $\rho_\alpha < \rho_0$, 
it would not be surprising that 
for the alternative distribution $\rho_{\tilde{\alpha}} > \rho_0$.
For a large range of initial conditions,
even if $\rho_\alpha < \rho_0$, 
the initially observed dynamics is rather the one 
described in Propositions \ref{p:r1r0ra} and  \ref{p2:r1r0ra}
by the interplay
between $1$, $0$ and $\tilde{\alpha}$ 
where $\rho_1 < \rho_0 < \rho_{\tilde{\alpha}}$.
The first step of the dynamics is thus 
obviously a convergence to 0
with a rate expected to stabilize towards $\rho_{\tilde{\alpha}} - \rho_0$.

\subsubsection{Validation by some simulations ? }
In order to evaluate the exceptionality 
of transitions towards 1,
one can propose the following simulation experiment,
which is a work in progress.
Let us consider some parameters 
for which most of the marginals tend towards 0,
with $\sigma$ sufficiently small.
To observe the upheaval of the group population
by type $G$ groups,
the marginals should be encoded rather not by the masses
at the different grid points, 
but by the logarithms of these quantities.
The case of the absorbing states is treated separately
from the marginal restricted to $(0,1)$.
For a better accuracy , 
it may be useful to refine the grid
in the vicinities of 0 and 1.

We do observe the mass towards 1
increasing up to the point of exceeding
the mass towards 0.
Yet, it is then unclear 
at which concentration the regeneration of these cooperative groups
exceeds the effect of having more and more exceptional transitions 
leading there.
Since the quickest transitions are expected 
to bring less mass towards $0$,
our idea is to truncate densities to prevent
the most exceptional transitions.
So at each simulation step,
we suppress from the marginal the mass 
on states than contain less than the threshold.
By varying the threshold,
we should have a better view 
on the number of groups required 
to observe the transition from $0$ to $1$
in individual-based models.

If the dynamics is almost unchanged after the truncation,
we could conclude that the cost of the transitions 
that mainly contribute to this upheaval 
is smaller than the threshold.
If the upheaval arises later on, 
it would mean that less costly transitions
 could have been sufficient to make type $G$ emerge.
Yet, their contribution becomes negligible 
when compared to quicker transitions.
Finally, 
if the marginal becomes supported 
on some interval that does not approach $1$,
it means that any transition towards $1$ 
would be at least as costly as the threshold.

\subsection{Conclusion of the section}
Generally, 
the interplay between different traits
happens in the time-scale 
at which their carriers 
can be differentiated.
Yet, we have seen in this example
that the trade-off between different kinds of advantages
can be particularly tricky.
The a priori neutral genetic drift 
might happen to be strongly 
coupled to the efficiency of some components of selection.
In a broader view, 
we can see this model as an illustration 
that selective effects
might be strongly dependent 
upon details of the local ecological dynamics,
and not only upon the average behavior.
If the local subdivision 
constitutes a sufficiently stable entity
with the ability to reproduce itself,
natural selection may act. 
Its strength
depends on the level of variability 
between those entities,
as if they were individuals.

Crucial requirements 
for the presented confrontation
are the independence of the groups 
bewteen successive splittings
and 
a strong heritability 
for the groups when they split
(the two descendants are assumed to be very similar 
to the "parent").
Even small interactions 
between these communities
(notably migration between groups)
is known to greatly disrupt the stability 
of cooperative strategies 
(cf. for instance \cite{WGGD06}).
So we clearly do not claim
that such simplified effect is prevalent,
because such lack of interaction
between the local dynamics is not so common.

Also, 
the selection at the group level
might rely on exceptional transitions
of the process $X$.
Although this exceptionality 
can be compensated in the long run
by a larger asymptotic growth rate,
two aspects should be remembered :
first, it might be much too unlikely
 for actual populations
 that some groups experience those rare transitions,
so that the mathematical model could be misleading;
second, it takes a very long time
for their descendants to invade.
It thus raises the question concerning real life
whether no other event happens to disturb
any of the sub-populations 
before the emergence of cooperative groups.

Finally, 
this competition model provides
a fruitful insight in an evolutionary perspective
(that is the subject 
of the following Section \ref{sec:tran}).
A main quantity of interest
is notably the probability	
that the descendants 
of a single mutant individual
with trait $y$ 
invade and replace the whole population 
of "residents" with trait $x$.
This probability is usually compared to the one for the invasion by mutants
identical to the residents (the neutral case).
In this model, for small values of $\sigma$,
as long as the mutant trait has a strong deleterious effect,
either at the individual- or at the group-level,
the invasion is much more difficult than neutral :
for the invasion of cooperatives, 
the random fluctuations
 inside the groups 
have to lead the process away 
from the very stable 
cheater quasi-equilibrium;
while, for the invasion of cheaters,
 the genetic drift between groups 
 has to disrupt 
 the also very stable cooperative quasi-equilibrium.
Such large deviations are known 
to be generally especially costly.
It means that natural selection of such traits
should be much more constraint.

\section{Large deviation estimate for the adaptation by mutations of weak effect}
\label{sec:henry}

In this section, 
I present the work of Champagnat and Henry
on some effects of Dirac concentration
in non-local models of adaptation with several resources \cite{Henry}.
The focus is on the trajectories of evolution,
with the selection of favorable mutations.

Interestingly, the authors of \cite{Henry}
show a mathematical similarity between two limiting behaviors :
a first one  where the trait variations are very small;
and the second
where the mutation rate is very small
(as compared to the selective effects).
Notably,
the first case 
can  be associated to a limit 
where mutations of very small effects accumulate.
The main requirement for the proofs
is that natural selection is very strong
in the timescale where one observes
the trait dispersion under neutrality.

Again, 
the population size is assumed
to be sufficiently large to include the whole range 
of the stochastic variations of the trait. 
The resulting purely deterministic model
shall provide a valid approximation 
to the dynamics of trait proportions 
in this population.
This relation to the trait proportions
can be retrieved from the convergence result
deduced from the time-scale separation.

\subsection{The continuous-space limiting behavior}
One considers the family $(u^\epsilon)_\epsilon$ 
of deterministic solutions  to the parabolic SDE :
 \begin{align}
&\partial_t u^\epsilon(t, x) 
 = \frac{\epsilon}{2}\, \Delta u^\epsilon(t, x) 
 + \frac{R(x, \psi_t^\epsilon)}{\epsilon}\, u^\epsilon(t, x),
 \quad x\in \mathbb{R}^d\, , \, t\ge 0,
 \label{uEps2}
 \\
&-\epsilon \log(u^\epsilon(0, x)) 
= h_\epsilon(x)
\notag
 \end{align}
where the competition effect $\psi_t^\epsilon
= (\psi_t^{i, \epsilon})_{i\le r}$ 
is defined, for $r$ resources,
with a competition kernel 
$\Psi_i: \mathbb{R}^d\rightarrow \mathbb{R}_+$ for resource $i$
 by  :
 \begin{align}
 \textstyle
 \psi_t^{i,\epsilon} := \int_{\mathbb{R}^d} \Psi^i(y)\, u^\epsilon(t, y)\, dy.
 \label{vEps}
 \end{align} 
 In this model,
 the traits $x\in \mathbb{R}^d$ characterize the ability to exploit the different resources.
Considering mutation effects through some heat kernel
corresponds to the case where
lots of mutations with very small effects
occur and are dispersed throughout the population.
In asexual populations, 
we rather expect 
selective mutations to invade and fix 
one after the other. 
Such a model is thus rather justified 
for sexual populations,
where many recombination events occur 
and many different alleles 
with small selective effects 
may coexist.
This is the principles 
of the so-called infinitesimal model
for which we refer notably to \cite{BEV17}.
%
For instance, 
it is known 
that many Human traits
are affected by many different loci in the genome.
Two individuals with almost the same phenotype
have possibly very different alleles in these loci
and the recombination of alleles 
along the lineages
creates variability.
The heat kernel seems then quite relevant to characterize 
the variability of response for such a trait.
It would be of interest to prove
that the dynamics of $u^\epsilon$,
i.e. as a solution to \ref{uEps2},
can be justified as a limiting description 
of individual-based models
as in Theorem \ref{MLim}.
 
The following assumptions are required in \cite{Henry} 
in order to justify a limit to 
$\varphi_\epsilon := -\epsilon \log u^{\epsilon}$ 
through a variational representation.
\\

\textsl{1. Assumptions on $\Psi_i$ : }
For any $1\le i\le r$, $\Psi_i \in W^{ 2,\infty}(\mathbb{R}^d)$.

Moreover, there exist $0< \Psi_{\min} <  \Psi_{\max}$ such that
$\frl{1\le i \le r}, \frlq{x\in \mathbb{R}^d}
\Psi_{\min} \le \Psi_i(x) \le  \Psi_{\max}$.
\\

\textsl{2. Assumptions on $R$ : }
(a) $R$ is continuous on $\mathbb{R}^d \times \mathbb{R}^r$.

\noindent (b) 
There exists $A>0$ such that
\begin{align*}
\frl{1\le i\le r} \frl{x\in \mathbb{R}^d}\frlq{(v_\ell) \in \mathbb{R}^r}
-A \le \partial_{v_i}R(x, (v_\ell)_{\ell \le r})\le -A^{-1}
\end{align*}

\noindent (c) 
There exist two positive constants $0< v_{\min} < v_{\max}$ such that
$\min_{x \in \mathbb{R}^d}
R(x, v) > 0$ as soon as $\|v\|_1 < v_{\min}$, 

and $\max_{x \in \mathbb{R}^d}
R(x, v) < 0$ as soon as $\|v\|_1 > v_{\max}$,
where $\|v\|_1 = \sum_{\ell\le r} |v_\ell|.$

\noindent (d) 
Let $H$ denotes the annulus $B (x, 2 v_{\max}) \setminus B(x, v_{\min}/2)$ (for the $\|.\|_1$ norm). Then
$\sup_{v\in H}
\|R(·, v)\|_{W^{2,\infty}} < \infty$.
\\
 
\textsl{3. Assumptions on $h_{\epsilon}$ :}
(a) $h_{\epsilon}$ is Lipschitz-continuous on $\mathbb{R}^d$,
 uniformly with respect to $\epsilon > 0$.
 
\noindent (b) $h_{\epsilon}$ converges 
uniformly
 as $\epsilon$ tends to 0 
 to a function $h$.
 
\noindent(c) For all  $\epsilon> 0$  and all $1 \le i \le d$,
$v_{\min} 
\leq \int_{\mathbb{R}^d} \Psi_i(x) \exp(- h_{\epsilon}(x)/\epsilon) dx
 \leq v_{\max}$.
 
In particular, $u_{\epsilon}(0, x)$ is bounded in $L^1(\mathbb{R}^d)$.
\\

\textbf{Remark : }
Such competition may at first sight
 appear as a global competition
where every trait is competing against each other for the same resources.
The first assumption done in \cite{Henry}
 requires indeed 
for the $\Psi_i$ to be lower-bounded by a strictly positive constant.
However, 
in a limit where $r$ is large and $\Psi_i$ is very concentrated,
the model may become a valid approximation of a local competition 
(in the trait space).
Thus, one may hope to extend the model naturally to such a framework.

\subsection{Results for the continuous case}
\label{sec:cont}
I gather here the main results presented in \cite{Henry},
where I have changed some notations for clarity,
notably $\psi_t$ and $\varphi$ :

\begin{theo}
\label{FKrep}
 ({\it Feynman-Kac representation of the solution of} (\ref{uEps2}))

 {\it Let} $u^{\varepsilon}$ {\it be the unique weak}
{\it solution of} (\ref{uEps2}), {\it then,
under the assumptions given in \cite{Henry} 
:}
\begin{align}
\ \forall(t, x)\in \mathbb{R}_+\times \mathbb{R}^{d},\quad 
\displaystyle
u^{\varepsilon} (t, x)
=\mathbb{E}_{x}\ \left[ \exp\left( \frac{-h_{\varepsilon}(X_{t}^{\varepsilon})}{\varepsilon}+\frac{1}{\varepsilon}
\int_{0}^{t}R(X_{t}^{\varepsilon}, \psi_{t-s}^{\varepsilon})ds
\right) \right],
\label{FKcont} 
\end{align}
{\it where for all} $x \in \mathbb{R}^{d}, \mathbb{E}_{x}$ {\it is the expectation
 associated to the probability measure} $\mathbb{P}_{x}$, {\it under which}
$X_{0}^{\varepsilon}=x$ {\it almost surely and the process} $B_{t}=(X_{t}^{\varepsilon}-x)/\sqrt{\varepsilon}$ {\it is a standard Brownian motion in} $\mathbb{R}^{d}.$
\end{theo}

\textbf{Remark :} 
This representation of the solution through a penalization
of a stochastic process 
is in practice of the same kind 
as the one presented for the empirical law $\mu_t$
in Section \ref{sec:aure}.
There is only a little difference
in the fact that $X_t$ describes now 
the process when we look at the lineage backward in time.
A similar convergence result 
should thus justify 
the interpretation of $u^\epsilon(t,x)\, dx$
as the limit of some individual-based measure-valued processes
$(\nu_t^{\epsilon, K})$.
In particular, 
$(1/\epsilon)\,.\, R(x, \psi_t^{\epsilon, K})$ 
should be the additional growth rate of individuals with trait $x$,
where for $i\le r$, 
$\quad \psi_t^{i, \epsilon, K} 
:= \langle \nu_t^{\epsilon, K}\, \vert \Psi^i\rangle$.
With fixed $\epsilon$, 
this is a priori  a specific case of the McKean-Vlasov equations 
mentioned in the conclusion of Section \ref{sec:aline},
where the law of the process itself
acts on the individual dynamics 
(again, I refer to \cite{S91}). 

\begin{lem}
\label{Phi}
{\it The function} $I^{R, \varepsilon}_t$ : $C([0, t])\rightarrow \mathbb{R}$
 {\it defined by}
$\quad
I^{R, \varepsilon}_t(y)=\int_{0}^{t}R(y_{s}, \psi_{s}^{\varepsilon})ds
\quad $
{\it is Lipschitz continuous on} $C([0, t])$ {\it endowed with the} $L^{\infty}$-norm.
 {\it The Lipschitz constant is uniform with respect to} $\varepsilon$ {\it for} $\varepsilon$ {\it small enough}.
{\it Moreover, there exists a kernel} $\mathcal{M}$ {\it on} $\mathbb{R}_+\times \mathcal{B}(\mathbb{R}^{k})$ {\it such that, along a subsequence} $(\varepsilon_{k})_{k\geq 1}$ {\it converging}
{\it to} $0$~:
\begin{center}
$\forall y\in C([0, t])
\quad 
I^{R}_t(y) 
:=\displaystyle \lim_{k\rightarrow\infty}  I^{R, \varepsilon_k}_t (y)
=\int_{0}^{t} \int_{\mathbb{R}^{k}} R(y_{s}, \psi)
 \mathcal{M}_{s}(d\psi) ds.$
\end{center}
\end{lem} 

By a kernel $\mathcal{M}$,
we specify  here a function from $\mathbb{R}_+\times \mathcal{B}(\mathbb{R}^{k})$ 
  into $\mathbb{R}_+$  such that, 
  for all $s\in \mathbb{R}_+, \mathcal{M}_{s}$ is
 a measure on $\mathcal{B}(\mathbb{R}^{k})$ 
 and, for all $A\in \mathcal{B}(\mathbb{R}^{k})$ ,  the function $s\rightarrow \mathcal{M}_{s}(A)$  is measurable.
\\

\begin{theo}
\label{Vtheo}
{\it For all} $(t, x)$ {\it in} $\mathbb{R}_{+}\times \mathbb{R}^{d},$
$$
\varphi(t, x)\ :=\lim_{k\rightarrow\infty}\ \varepsilon_{k}\log u^{\varepsilon_{k}}(t, x)
=\ \underset{y\in \mathcal{G}_{t,x}}{\sup}\ \{-h(y_{0})+I^R_t (y)-I^L_{t}(y)\}
$$
{\it where the convergence holds uniformly on compact sets and the limit}
$\varphi(t, x)$ {\it is Lipschitz} $w.r.t.$\ $(t, x)\in\mathbb{R}_+\times \mathbb{R}^{d}$
{\it while}\quad 
$I^L_{t}(y)=\left\{\begin{array}{ll}
\frac{1}{2}
\int_{0}^{t}\Vert y'(s)\Vert^{2}ds & if\ y\ is\ absolutely\ continuous,\\
+\infty & otherwise.
\end{array}\right. $

$\mathcal{G}_{t,x}$ {\it denotes the set of continuous functions from} $[0, t]$ {\it to} $\mathbb{R}^{d}$ {\it such that} $y_{t}=x$,
 {\it and} $I^R_t$ 
 {\it and} $(\varepsilon_{k})_{k\geq 1}$ {\it are associated as in Lemma \ref{Phi}.} 
\end{theo} 

\textbf{Remarks :}
$\bullet$ In this optimization,
$y$ defines the main ancestral line of the individuals
dominant with trait $x$ at time $t$.
Namely, 
the Large Deviation theory ensures 
that as $\epsilon \rightarrow 0$, 
their ancestral lines are very concentrated 
on such specific and deterministic histories.

$\bullet$  
Given a lineage described by $y$, 
$I^R(y)$ encodes an average effect of selection 
due to the perceived growth rate.
As one can infer from (\ref{FKcont}) when $\epsilon$ is very small,
the density around such a path is approximately amplified
by the exponential of this quantity divided by $\epsilon$.
Such large amplification is the core 
of the Large Deviation theory, 
where the stochastic behavior 
is concentrated close to specific paths.

$\bullet$  
The dependency of the competition kernel
$\mathcal{M}_{s}$ 
on $\varphi$ is implicit and unclear 
in the limiting model.
The level of competition depends 
on the global composition 
in traits of the population,
while $\varphi$ only indicates which traits are non-negligible 
(the $x$ for which $\varphi(x) = 0$).
Detailed information on the density
is lost in the limit $\epsilon \rightarrow 0$,
and that is why subsequences $\epsilon_k$ 
are considered. 
For the same reason, 
the uniqueness of solutions to such equations 
is not easy to establish
(when it holds).
\\

Again in \cite{Henry},
 the authors finally relate
this variational limit 
to the solution of the Hamilton-Jacobi equation~:
\begin{theo}
\label{HJtheo}
Under some assumptions on the behavior of $R$  and  $h_\epsilon$,
$\psi^{\varepsilon_{k}}$ {\it converges in} $L_{loc}^{1}(\mathbb{R}_{+})$
{\it along the subsequence} $\varepsilon_{k}$ {\it of Lemma} \ref{Phi} {\it to a nondecreasing limit} $\overline{\psi}$,
{\it the kernel} $\mathcal{M}$ {\it satisfies}
\begin{center}
$\forall s\geq 0,\qquad   
\mathcal{M}_{s}(d\psi)=\delta_{\overline{\psi}_s}(d\psi),$
\end{center}
{\it and the limit} $\varphi$ {\it of Theorem \ref{Vtheo}} {\it solves in the viscosity sense~:}
\begin{center}
$\left\{\begin{array}{ll}
\partial_{t}\varphi(t, x) =R(x,\overline{\psi}_{t}) 
+\frac{1}{2}|\nabla \varphi(t, x)|\ ,
 \quad& \forall t\geq 0, x\in \mathbb{R}^{d},\\
\max_{x\in \mathbb{R}^{d}}\varphi(t, x)=0,\ &\forall t\geq 0.  
\end{array}\right.$
\end{center}
\end{theo}

\paragraph{\textsl{Behavior of the solution to Hamilton-Jacobi equation}}
\textcolor{white}{.}

From Theorem \ref{Vtheo}, 
one deduces the following approximation
for small $\epsilon$ :
$u^\epsilon(t, x) \approx C(t, x)\, \exp[\varphi(t,x)/\epsilon]$.
We thus expect to observe
a concentration of the individual traits 
on typically a unique value 
or a few well-separated values.
In the limiting model, 
these values may be driven by a continuous displacement
in the direction of natural selection.
They may also split, 
which would in a sense describe
 a model of speciation.
It is in fact the original purpose of such Hamilton-Jacobi analysis
in \cite{DJMP05}
to extend the model of Adaptive Dynamics (cf. Subsection \ref{AD}).
The reverse can happen too,
with two subpopulations 
merging into a single one :
 the traits are 
 concentrated around a moving value
 specific to each subpopulation
 until these two values coalesce.
Finally, there could also be some jumps, 
that is a brutal change of the trait composition.
A wide variety of events are thus observable 
with this simplified model.
\\

\paragraph{\textsl{Discussion on the involved time-scales}}\textcolor{white}{.}

At a given $\epsilon>0$,
looking at equation (\ref{uEps2}),
we can infer from the growth rate 
that the stabilization of the densities
occurs in a time-scale of order $\epsilon$,
as compared to the time-scale at which the traits evolve.
The variance of the trait evolves on the contrary 
at a much larger time scale of order $\epsilon^{-1}$,
as in the case of neutral traits.
Of course, the detailed description 
of the demography
as presented in 
Section 1 above and in \cite{Aline}
should occur in the short time-scale of order $\epsilon$.
In such a large population, 
the life-expectancy of any particular individual
is on an even shorter time-scale.
depending on the ratio between the birth and the death rates.
The closer it is to one, the shorter this time-scale.
A separation of time-scale with the evolutionary trajectory
can still be justified 
for the competition between different effects of selection.
But as we have seen in Section \ref{sec:aure},
this competition might last for very long times.
\\

In practice, 
the description given by equation (\ref{uEps2})
is not so easily related to 
the macroscopic observation 
of the response to selection,
because the neutral variability is difficult to scale
(the factor $\epsilon$ before the Laplacian).
We can reasonably assume 
that we know the density $ u^e(t,x) = u^\epsilon(t,x)$ 
(as long as it is not too negligible)
and the fitness effect $R^e_t(x) = (1/\epsilon)\,.\, R(x, \psi_t)$.
For at least positive values,
it can be inferred in some laboratory experiment
by artificially introducing 
individuals with trait $x$
in a much larger sample of the population at time $t$
and see how they grow.
One could then estimate the value $\epsilon$
for which $\varphi^\epsilon(t,x) = \epsilon\log u^e(t,x)$
seems reasonably close to the solution $\varphi$
of the Hamilton-Jacobi equation :
\begin{equation}
\textstyle
\partial_{t}\varphi(t, x) 
= \epsilon\, R^e_t(x) 
+\frac{1}{2}|\nabla \varphi(t, x)|,
\label{HJe}
\end{equation}
or the associated control problem.

Note that a Gaussian distribution for $u^e$ 
 corresponds to $\varphi(t, x) \sim - (x- x_t)^2/\sigma_t^2$,
 that is the Second order Taylor expansion around $x_t$,
 with $x_t$ the optimal trait 
 and $\epsilon\, \sigma_t^2$ the variance of the distribution.
Moreover, if $\varphi$ is a solution to equation $(\ref{HJe})$,
then for $\lambda >0$, 
$\hat{\varphi}(x, t) := \varphi(\lambda\, t, \lambda\, x)$
is a solution to $(\ref{HJe})$ with $\epsilon$ 
replaced by $\lambda\, \epsilon$.
Since we expect the dependency in $\sigma_0$ on $\sigma_t$
to vanish very quickly,
the identification of $\epsilon$ from $\sigma_t$
should not be a too difficult option.
Yet, as mentioned in the next paragraph,
the speed of response to selection 
might be driven by an exceptional proportion 
of the population at time $t$.
An estimate such as $\sigma_t$
that summarizes rather the core of the distribution $u^e(t,x)\, dx$
might thus not be so relevant.
\\

\paragraph{\textsl{Connection with individual-based models}}\textcolor{white}{.}

To justify the connection with individual-based models,
concentration effects may appear favorable.
Yet,
as we can see from Theorem \ref{Vtheo},
the adaptation of the population is driven 
by exceptional profiles of stochastic variations.
The stochastic techniques of Large Deviations 
is then crucial to obtain the limiting behavior.
Notably, it means 
that selection is mainly driven by individuals at the front.
This optimization provides a very interesting insight
on the history of the genealogies 
that led to the traits observed generally in the population 
at time $t$.
Notably,
we can  observe cases 
where the  traits at time $s\le t$ for the ancestors of the dominants at time $t$
may not be typical at all as compared to the dominant traits at time $s$.
This of course raises the issue of 
the biological relevance
 of such transitions.
We face similar limitations  
as in previous Section \ref{sec:aure}
for the relevance of events 
involving very small population sizes.
Still, 
it is much easier to relate 
the above description
to some individual-based model
than from the more classical description 
only relying on Hamilton-Jacobi equations
(cf. the conclusion of the section).

\subsection{Results for the discrete case}
\label{sec:disc}
This analysis of Feynman-Kac operators and Large Deviation principles
can be adapted to a discrete space.
Simplifying assumptions are then more easily obtained
and makes it possible to have a clearer view on the implications of this model.

In \cite{Henry},
the authors consider
the following system of ordinary differential equations~:

$\left\{\begin{array}{ll}
&\displaystyle
\partial_t u^{\varepsilon}(t, k)
=\sum_{j\in E\backslash \{k\}}
\exp\big( \frac{-\mathfrak{T}(k,j)}{\epsilon}\big)
(u^{\varepsilon}(t, j)-u^{\varepsilon}(t, k))
+\frac{1}{\varepsilon} u^{\varepsilon}(t, k) R(k, \psi_{t}^{\varepsilon})\ ,
 \hspace{0.3cm} \forall t\in[0, T]\ ,\ \forall k\in E
 \\&
u^{\varepsilon}(0, k)
=\exp(-\frac{h(k)}
{\varepsilon})\ , 
\qquad \psi_{t}^{i, \varepsilon}
= \sum_{k\in E} \Psi^i(k)\, u^{\varepsilon}(t, k),\quad \forall i\le r
\end{array}\right.$
\\
In such a model, 
mutations from state $j$ to state $i$ 
happen at rate $\exp(-\epsilon^{-1} \mathfrak{T}(k,j))$.
As mentioned in the introduction of \cite{Henry},
this is not the classical form 
of the solution for the system of ODE 
describing population densities.
The former can however be obtained
from the latter
by a slight adjustment in the definition of $R$
(vanishing as $\epsilon \rightarrow 0$)
and is more practical for the following analysis.
The growth rate of the individuals
in state $k$
is (nearly) $\epsilon^{-1} R(k, \psi)$, 
where $\psi = (\psi^i)$ 
specifies the amount of available resources.


First,
the solution $u^{\varepsilon}$ of the system
 is described in \cite{Henry}
by using an integral representation similar to  (\ref{FKcont}). 
Let $(X_{s}^{\varepsilon},\ s\in[0, T])$ be the Markov process in $E$ with infinitesimal generator~:
$$\textstyle
L^{\varepsilon}f(k)=
\sum_{j\in E}
\exp[\frac{-\mathfrak{T}(k,j)}{\varepsilon}]
(f(j)-f(k))
$$
i.e. the continuous-time Markov process 
whose jump rate from state $i\in E$ 
to $j\neq i$ is
 $\exp(-\mathfrak{T}(i,j)/\varepsilon)$ .

\begin{prop}
({\it Integral representation}) 

{\it For any positive real number} $t$ {\it and any element} $i$ {\it of} $E,$
{\it we have}
$$
u^{\varepsilon}(t, i)
=\mathbb{E}_{i}\ \left[ \exp\left( \frac{-h(X_{t}^{\varepsilon})}{\varepsilon}
 +\frac{1}{\varepsilon}\int_{0}^{t} R\ (X_{s}^{\varepsilon}, \psi_{t-s}^{\varepsilon})ds
 \right)\right]
$$
\end{prop} 
The interpretation is  very similar to the one of Theorem \ref{FKrep}.

\begin{prop}
({\it Weak LDP})
 $(X^{\varepsilon})_{\varepsilon\geq 0}$ {\it satisfies a weak LDP with rate function}
\begin{center}
$I_{T}:\ 
\left\{\begin{array}{ll}
&\mathrm{D}([0, T], E)\ \rightarrow\ \mathbb{R}
\\
&y\ \mapsto\ \sum_{\ell=1}^{N_{y}}\mathfrak{T}(y_{t_{\ell}^{y}-}, y_{t_{\ell}^{y}})\ ,
\end{array}\right.$
\end{center}
\noindent
{\it where} $\mathrm{D} ([0, T], E)$ {\it is the space of c\`{a}dl\`{a}g functions from} $[0, T]$ {\it to} $E$,
 $N_{y}$ {\it is the number of jumps of} $y$
 {\it and} $(t_{\ell}^{y})_{1\leq \ell\leq N_{y}}$ 
 {\it the increasing sequence of jump times of} $y.$
\end{prop}

We shall use the notation
$\quad
I_{T}(y)
:=\sum_{0<s\leq T}\mathfrak{T}(y_{s-}, y_{s})
\quad$
with the implicit convention that $\mathfrak{T}(i, i) =0$ for all $i\in E.$

\begin{theo}
\label{Dtheo}
Lemma \ref{Phi} is also satisfied in the discrete case and
{\it for all} $(t, i)$ {\it in} $(0, +\infty) \times E,$
with the associated subsequence $(\varepsilon_{k})$~:
\begin{align*}
&\varphi(t, i)\ :=\lim_{k\rightarrow\infty}\varepsilon_{k}\log u^{\varepsilon_{k}}(t, i)
\\&
=\underset{y\in \mathrm{D}([0,t],E) s.t.\ y_0=i}{\sup}\{-h(y_{t})+\int_{0}^{t}\int_{\mathbb{R}^r}
R(y_{s},\psi)\mathcal{M}_{t-s}(d\psi)ds-\sum_{0<s\leq t}\mathfrak{T}(y_{s-}, y_{s})\}\ .
\end{align*}
\end{theo} 

As in the expression in Theorem \ref{Vtheo},
one shall remark 
that the optimal $y$ in this expression
describes the transition of states 
for the lineage of a typical individual at time $t$.
This evolution is backward in time,
as can be also seen 
in the term $R(y_{s},\psi)\mathcal{M}_{t-s}(d\psi)$,
where $\mathcal{M}_{t-s}(d\psi)$ provides the law of $\psi_{t-s}$.

\begin{theo}
{\it For all subsequence} $(\varepsilon_{k})_{k\geq 1}$ 
{\it as in Theorem} \ref{Dtheo},
 {\it the limit} $\varphi(t, i)$ 
 {\it of} $\varepsilon\log u^{\varepsilon}(t, i)$ {\it is}
{\it Lipschitz with respect to the time variable} $t$ {\it on} $(0, +\infty)$ . {\it In addition, if} $h(i) \leq h(j)+\mathfrak{T}(i, j)$ {\it for all}
$i\neq j$, {\it the function} $\varphi$ {\it is Lipschitz on} $\mathbb{R}_{+}.$
\end{theo} 

We see that all the previous results in continuous space 
have an equivalent for discrete space.
\\

In the following, 
we focus on a specific case where the uniqueness can be obtained
and the whole dynamics described with much more accuracy.
The result presented in \cite{Henry}
relies on some assumption on the stability of the dynamics
restricted to any subset of $E$, 
which is called Assumption $(H)$. 
I present next the main intuition behind this assumption
and refer to \cite{Henry} 
for the exact formulation.

By this Assumption $(H)$, 
the authors notably ensure that,
for any subset $A$ 
containing $k$ different types,
there exists a unique strongly attractive equilibrium.
Since the focus is on the dynamical system, 
the estimation of the dynamics is valid 
as soon as the $k$ types are initially in non-negligible proportion.
Yet, 
some of the components of the steady-state may be 0,
meaning an exponential decay of the mass of those.
Assumption $(H)$ also ensures that any other (unstable) equilibrium 
is quickly escaped.
The dominant traits at time $t$ (for which $\varphi(t, .) = 0$)
are then proved to stay piecewise constant,
while the emergence of favorable competitors is easily computed.
Once such a competitor has emerged,
it disrupts the equilibrium of traits,
and assumption $(H)$ enables us to predict the issue
of the equilibrium that follows.

Given any set $A$ of present types,
the equilibrium is given by $(u^*_{A, j})_{j\in A}$.
Thus, the competition exerted on the $i$-th resource by this eco-system
is $\quad \Psi^i(A) := \sum_{j=1}^{r}\Psi^{i}(j)\,.\, u_{A,j}^{*}.$
\\

\textsl{Remark : }
To satisfy assumption $(H)$,
these steady states need to satisfy
a condition of compatibility 
regarding their vanishing components.
By restricting the dynamics 
on some subset $B$ of $A$
that contains 
all the non-vanishing components 
of the steady state associated with $A$,
the steady state for $B$
is necessarily the restriction on $B$ 
of the steady state for $A$.

\begin{prop}
 {\it Assume that} ({\it H}) {\it is satisfied}. {\it Let} $(\varepsilon_{k})_{k\geq 1}$
  {\it be as in Theorem} \ref{Dtheo}. {\it For any} $t\geq 0$, {\it there}
{\it exists} $\rho_{t}>0$ {\it such that, 
for all} $ s\in (t, t+\rho_{t}$] ,
 $\psi_{s}^{\varepsilon}$ {\it converges to} 
$\psi_s = \psi_t := \Psi(\{\varphi(t, \cdot)=0\})$,
{\it where the convergence}
{\it is uniform in all compact subsets of} $(t, t+\rho_{t}]$
 {\it and where we define for any $A\subset E$~:}
$\quad F(A)
=\ (\sum_{j=1}^{r}\eta_{i}(j)u_{A,j}^{*})_{1\leq i\leq r}.$

{\it In particular, the weak limit} $\mathcal{M}_{s}$ {\it of} $\delta_{\psi^{\epsilon_k}_s}$ 
{\it obtained in Lemma} \ref{Phi} 
{\it satisfies}
$\mathcal{M}_{s}=\delta_{\psi_t}$, 
{\it for almost all} $s\in(t, t+\rho_{t})$
{\it and the function} $t\mapsto \psi_t$ {\it is right-continuous}.
\end{prop}

As in the previous case,
we shall observe a concentration effect. 
It means here that 
there are brutal transitions 
between one demographic equilibrium
and the next one generated 
by the first invasion by a mutant sub-population.
\\

\begin{prop}
 {\it Assume that} ({\it H}) {\it is satisfied}. 
 {\it Any limit} $\varphi$ 
 {\it of} $\varepsilon_{k}\log u^{\varepsilon_{k}}$ 
 {\it along a subsequence as in Theorem} \ref{Dtheo}  
 satisfies {\it for all} $i\in E$ 
 {\it and for all} $t\geq 0$ : 
 $\varphi(0, i)=-h(i)$ 
$$
\varphi(t, i)=\sup_{y(0)=i}\{-h(y(t))+\int_0^{t} R(y(u), \psi_{t-u})\ du 
-I_{t}(y)\},
$$
{\it and its dynamic programming version}
\begin{align}
\varphi(t,\displaystyle \ i)
=\sup_{y(s)=i}\{\varphi(s, y(t))+\int_0^{t} R(y(u), \psi_{t+s-u})\ du-I_{s,t}(y)\}.
\label{varVf}
\end{align}
{\it In addition, the problem (\ref{varVf}) admits a unique solution s.t.} $t \mapsto \psi_t := F (\{\varphi(t, \cdot) = 0\})$ {\it is right-continuous. 
In particular, the full sequence} $(\varepsilon\log u^{\varepsilon})_{\varepsilon>0}$ {\it converges to this unique solution when} $\varepsilon\rightarrow 0.$
\end{prop} 

\begin{theo}
 {\it Under Hypothesis} ({\it H}) {\it and assuming that} $h(i) \leq h(j) +\mathfrak{T}(i,j)$ {\it for all} $i \neq j$, {\it the}
{\it problem}
\begin{align*}
\left\{\begin{array}{ll}
& \partial_t \varphi(t, i)=\sup\{R(j, \psi_t)\,
|\, j\in E\, s.t.\, \varphi(t,j)-\mathfrak{T}(j, i)=\varphi(t, i)\},
 \\
&\forall i\in E\, , \, \quad  \varphi(0, i)=h(i)
\hcm{2}(\text{ with the convention } \mathfrak{T}(i, i) = 0)
  \end{array} \right.
\end{align*}
 {\it admits a unique solution such that} $t\mapsto \psi_t = F(\{\varphi(t, \cdot) = 0\})$ {\it is right-continuous and it is the unique solution to the variational problem}  (\ref{varVf}).
\end{theo}

The system follows a succession 
of equilibria (with possibly different types that equilibrate),
with the population headcount of mal-adapted types
vanishing linearly in their logarithm,
while this headcount logarithm increases 
for the rare adapted ones.
It goes on 
until one adapted type 
reaches the threshold for a non-negligible frequency.
Then occurs a kind of "catastrophe", 
where the whole equilibrium is immediately renewed 
(in the time-scale of evolution).
Given the similarity 
between the description in this discrete case 
and the continuous case presented in Subsection \ref{sec:cont},
one can infer that similar "catastrophes"
might happen even in the continuous case.
This is less expected,
since we see meanwhile a dynamics for the dominant trait,
but it cannot be excluded a priori
that a brutal change of the traits happen.
\\

Without this assumption of a unique stable equilibrium, 
it would not be clear what happens to the population 
during such  a "catastrophe" event.
Namely, 
the only information provided by $\varphi(., t)$ 
is the knowledge of the non-negligible types, 
for which $\varphi(x, t) = 0$,
with a priori no means to infer otherwise the headcounts at equilibrium.
Notably, different issues or an unstable behavior
 could lead to different competition effects,
 thus different dynamics after the "catastrophe".
\\

The main restriction for this model is that
 the tails of the distribution are not too much involved.
Indeed, 
when there is some barrier,
 where the growth rate is very low,
the trajectory that gets selected upon 
through the optimization procedure in Theorem \ref{Vtheo}
might not be realistic.
It might well be inferred from this optimization
that  the ancestors of the population at time $t$
shall represent at time $s\le t$ 
a proportion less than $10^{-9}$ of the entire population at time $s$.
Of course, 
it would only be reasonable for populations 
with a size quite larger than $10^{9}$ !
And all the more since families 
with small population size
have a high risk of going extinct.
Considering the actual population size considered,
such optimization procedure could be an efficient way 
to truncate too rare such transitions
as we have proposed for the analysis in Section \ref{sec:aure}.
Namely, 
for a threshold $- \varphi_\vee$,
we would like to define some
 $\widetilde{\varphi}(t,x)$
as the supremum over the same quantity 
as in Theorem \ref{Vtheo}
under the additional condition on $y\in \mathcal{G}_{t,x}$
that :
\begin{align*}
\frlq{s\le t}
-h(y_{0})+I^R_s(y)-I_{s}(y) > - \varphi_\vee,
\quad \text{ where }	
I^R_s(y) := \int_{0}^{s} \int_{\mathbb{R}^{k}} R(y_{u}, \psi) \widetilde{\mathcal{M}}_{u}(d\psi) du.
\end{align*}
Naturally, we want to impose
$\widetilde{\varphi}(t,x) = -\infty$ 
if no function  $y\in \mathcal{G}_{t,x}$ 
is able to satisfy this condition.
By this way, we forbid transitions
that rely on too exceptional ancestors
to be treated in such a deterministic way.
While $\widetilde{\mathcal{M}}_u$
shall represent the availability of resources
at time $u$,
it shall actually 
depend on $\widetilde{\varphi}$,
as in the case without truncation.
This makes the analysis a priori much trickier.
Yet, the transitions we forbid are rare 
so there is a lag between the time 
at which they occur 
and the time 
at which they have an effect 
on $\widetilde{\mathcal{M}}$.
The estimation of $\widetilde{\mathcal{M}}$
 shall thus not be more difficult
than the one of $\mathcal{M}$
 in the case without truncation.

Alternative strategies have been proposed
 in \cite{PG09} and \cite{MBPS}
 with more regularities than the proposed truncation,
 yet not much more biological justification.
To prevent such exceptional densities,
 the authors use a singular term for the growth rate,
 which gets very negative when the density of the state is too small.

To be even more realistic, 
one should treat such events of crossing barriers 
as punctual events occurring at a very low rate.
In fact, 
large deviation approaches may be well 																														suited to estimate these rates.
One should only remark that 
the growth rate on the other side of the barrier
is not responsible for any increase
on the rate at which such transitions occur.
It only increases the probability 
that such exceptional crossing event 
leads actually to an invasion.
And that is why such events are much more exceptional 
than with the asymptotic model given by $\varphi$,
where the growth rate 
immediately and regularly
 increases the log-density on this other side.
The punctual transitions 
that we are to describe in the next section
is thus in fact well suited 
to deal also with such rare crossing events.

\section{Mutations as the limiting factor}
\label{sec:tran}

Although this section is clearly 
referring to the talk given by Tran,
the presentation he made was not dedicated
to a specific paper.
It was rather
a general introduction to the framework of IBM 
at use to prove
 classical scaling limits of adaptive dynamics.
 The focus was mostly on models 
 where the traits are highly conserved along the lineages
 until exceptional events of mutation
 in a new-born.
Relying on two articles of my choice 
that are in the spirit 
of Tran's talk,
I intend to offer a broad perspective
 on the limiting descriptions of IBM.
Given the progression of the whole paper,
I favored in my choice
(besides the commitment of the speaker)
the coupling of several time-scales 
over the intrinsic complexity in the asymptotic dynamics.
\\

This last Section 4 addresses the longest time-scales,
where the limiting factor for evolution 
is the emergence of mutations.
By its emergence,
we mean not only that some individuals
carry such mutations at time $t$
but rather that it has been stabilized at a non-negligible frequency
in the population.
For simplicity, 
many authors assume that selective traits cannot coexist
in the same population,
and that the population size is very stable 
at a value depending on this dominant trait.
Namely,
demographic fluctuations 
that has been crucial in Section \ref{sec:aure}
are neglected,
and a fortiori the individual variability detailed in Section \ref{sec:aline}
is only seen through an average effect.
In such a deterministic approximation,
the fixation or the decay of an invasive trait 
is generally described,
as soon as the proportion of invaders is not negligible anymore,
by Lotka-Volterra equations.
For elementary models of interactions 
(notably not frequency-dependent),
it is then easy to impose conditions ensuring 
that the fixation
of one of the alleles
 is the only stable equilibrium
of the dynamical system.

Given the assumption
that such ecological transitions 
happen in a much shorter time-scale
than the time between
occurrences of new mutations 
(at first in a unique individual),
such conditions, 
called "invasion implies fixation",
ensures that the population stays almost always
monomorphic in the evolutionary time-scale.
This means in particular
that the invasion of new mutations can be treated as punctual events
that refreshes the equilibrium.
In this view, 
this behavior is similar to
the evolution of traits given for the discrete space
in Section \ref{sec:henry}
(cf. the concluding remarks for more details on this link).
Yet, 
these punctual fixation events happen at unpredictable times,
corresponding to the occurrence of a successful mutation.

Note that this shall hold true
 even if coexistence of traits are allowed in the model.
Yet, the rate at which 
an invasion occurs
and its issue
might depend on any variations 
of the frequencies (like oscillation of frequencies)
and of the population size,
which makes the analysis tricky.
Of course, one could also generalize these models
to include competition
for resources as in Section \ref{sec:henry}.
Yet, 
the specificity of these more complicated models
can only be seen when several traits can coexist.
If one wishes to include events of speciation,
this is however a very reasonable way to justify it 
(specialization of two sub-populations).

For strong selection effects,
one can prove negligible 
the time-interval
between the arrival of the first mutant
and the invasion of the population
by its descendants
(at least reaching a non-negligible frequency).
Namely,
a favorable mutant has a larger expectancy of offspring
than individuals of resident type and shall thus generate a family 
that either dies out quite early
or thrive at an exponential rate.
Contrary to the model of discrete trait evolution
described in Subsection \ref{sec:disc},
the time-scale for the occurrence of mutations
is chosen here in such a way
that this growth period
is not a limiting factor.
By "successful mutation",
 we mean the occurrence of the mutation
in a first carrier
 followed by a successful invasion 
 by its descendants 
 of the resident population.
 By assuming that mutation events are sufficiently rare,
 we infer 
 that successful mutations 
 happen nearly independently of the time waiting for it 
(where the monomorphic population stays the same)
and of unsuccessful invasions.
This explains why the time-interval between two
successive events of invasion 
is given by an exponential law without memory.

Such a process with piecewise constant population equilibria,
including a dominant trait,
has been originally introduced 
as the Trait Substitution Sequence in \cite{M96}
and more formally related to individual-based models 
in \cite{C06}.
A large family of extensions has emerged,
notably to include coexistence  of traits 
(cf. \cite{CM11}, 
\cite{BB18}),
or more complicated interactions, for instance with horizontal transfer
\cite{BCFMT} 
or aging \cite{MT09}.
Since they deal with two additional time-scales,
we rather focus on \cite{bfmt18},
where events of invasion affect
 the dynamics of some marker,
 and \cite{BBC17},
 which demonstrates that adaptive dynamics
 can indeed be 
an accurate approximation of some individual-based models.

\subsection{A model combining a marker dynamics with the evolution of traits}

The next four subsections present 
the results of Billiard, Ferri\`{e}re, M\'el\'eard and Tran
on the stochastic dynamics of neutral markers 
coupled to the one of adaptive traits \cite{bfmt18}.
The authors  consider an asexual population 
driven by births and deaths 
where each individual is characterized by hereditary types: 
a phenotypic trait under selection and a neutral marker.
This analysis is motivated 
by research on the prevalence of selection
between various species or clades 
that could be based on the observed variability of neutral markers.

The trait and marker spaces $\mathcal{X}$ and $\mathcal{U}$ 
are assumed to be compact subsets of $\mathbb{R}$. 
The type of individual $i$ is thus a pair $(x_{i}, u_{i})$,
 $x_{i} \in \mathcal{X}$ being the trait value 
 and $u_{i} \in \mathcal{U}$ its neutral marker. 
The individual-based microscopic model 
from which we start 
is a stochastic birth and death process 
with density-dependence  
whose demographic parameters are functions of the trait under selection 
and are independent of the marker. 
We assume that the population size scales 
with an integer parameter $K$ 
tending to infinity 
so that individuals are weighted
with $\textstyle \frac{1}{K}$
to observe a non-trivial limit
of the empirical measure. 
The state of the population at time $t \geq 0,$
rescaled by $K$, is described by the point measure :

\begin{align*}
&\displaystyle \nu_{t}^{K}
=\ \frac{1}{K}\sum_{i=1}^{N_{t}^{K}}\delta_{(x_{i},u_{i})}
 =X_{t}^{K}(dx) \, \pi_{t}^{K}(x, du),
\quad  \text{where } 
 X_{t}^{K}=\frac{1}{K}\sum_{i=1}^{N_{t}^{K}}\delta_{x_{i}}
\text{ and }
\displaystyle \pi_{t}^{K}(x, du)
=\ \frac{\sum_{i=1}^{N_{t}^{K}}
\mathbf{1}_{x_{i}=x}\delta_{u_{i}}}
{\sum_{i=1}^{N_{t}^{K}} \mathbf{1}_{x_{i}=x}}
\end{align*}
are respectively the trait marginal
and the marker distribution for a given trait value $x$.
Here, $\delta_{(x,u)}$, $\delta_x$ are respectively the Dirac measure at $(x, u)$ and $x$. 

With the following definitions,
the authors ensure  that the mutations happen 
at different time scales for the trait and for the marker, 
both longer than the individuals lifetime scale. 
Thus, the limiting behavior results 
from the interplay of these three time scales:
 births and deaths, trait mutations and marker mutations.
To justify such a separation 
of time scales for the mutations,
the proof relies strongly on the fact
 that the population size remains
around the equilibrium of some dominant trait(s)
when a mutation on the trait occurs.
Although the authors are able to
include fluctuations of the parameter $u$,
it is only possible because $u$ 
has no effect 
 on the probability that the mutation succeeds to invade the population.

\subsection{The individual-based model}
\label{sec:iBm}

\begin{Def}
\label{DefT}
$\bullet$ {\it An individual with trait} $x$ {\it and marker} $u$ {\it reproduces with birth rate}
{\it given by} $0\leq b(x) \leq\overline{b}$, {\it the function} $b$ {\it being continuous
and $\overline{b}$ a positive real.}

$\bullet$ 
{\it Reproduction produces a single offspring 
which usually inherits the trait and marker}
{\it of its ancestor except when a mutation occurs. 
Mutations on trait and marker occur}
{\it independently with probabilities} $p_{K}$ 
{\it and} $q_{K}$ {\it respectively. 
Mutations are rare 
and the marker mutates much more often than the trait. 
We assume that}
\begin{center}
$q_{K}=p_{K}r_{K}$ , 
{\it with} $p_{K}=\displaystyle \frac{1}{K^{2}}, q_{K}\underset{K\longrightarrow\infty}{\rightarrow} 0$, 
$ r_{K}\underset{K\longrightarrow\infty}{\rightarrow}+\infty$.   
\end{center}
$\bullet$ {\it When a trait mutation occurs, 
the new trait of the descendant is} $x+k \in \mathcal{X}$ 
{\it with} $k$ {\it chosen according to the probability measure} $m(x, k)dk.$

$\bullet$ {\it When a marker mutation occurs, 
the new marker of the descendant is} $u+h\in \mathcal{U}$
 {\it with} $h$ {\it chosen 
 according to the probability measure} $G_{K}(u, dh)$.
{\it For any} $u\in \mathcal{U}, G_{K}(u, .)$ 
{\it is approximated as follows when} $K$ {\it tends to infinity}:
\begin{center}
$\displaystyle \lim_{K\rightarrow+\infty}\sup_{u\in \mathcal{U}}
\left |\frac{r_{K}}{K} \int_{\mathcal{U}} (\phi(u+h)-\phi(u)) G_{K}(u, dh)
-A\phi\right|\ =0$,  
\end{center}
{\it where} $(A, \mathcal{D}(A))$ 
{\it is the generator of a Feller semigroup and} 
$\phi \in \mathcal{D}(A) \subseteq 
C_{b}(\mathcal{U}, \mathbb{R})$,
{\it the set of continuous bounded real functions on} $\mathcal{U}.$

$\bullet$ {\it An individual with trait} $x$ 
{\it and marker} $u$ 
{\it dies with intrinsic death rate} $0\leq d(x) \leq\overline{d},$
{\it the function} $d$ {\it being continuous
and $\overline{d}$ a positive real. 
Moreover,
 the individual experiences competition the}
{\it effect of which is an additional death rate} 
$$
\eta(x)\,.\, C*\nu_{t}^{K}(x) 
= \displaystyle \frac{\eta(x)}{K}\sum_{i=1}^{N_{t}^{K}}C(x-x_{i}).
$$

{\it The quantity} $C(x-x_{i})$ {\it describes the competition pressure exerted by an individual}
{\it with trait} $x_{i}$ {\it on an individual with trait} $x$. {\it We assume that the functions} $C$ {\it and} $\eta$
{\it are continuous and that there exists} $\underline{\eta}>0$ {\it such that}
\begin{center}
$\forall x, y\in \mathcal{X},\ \eta(x)\ C(x-y)\ \geq\underline{\eta}>0$.   (2.6)
\end{center}
A classical choice of competition function $C$ is $C\equiv 1$ which is called 
"mean field case" or "logistic case". In that case the competition death rate is $\eta(x)N_{t}^{K}/K.$
\end{Def}

As mentioned in the beginning of the Section 4, 
the authors work under the simplifying assumption that follows,
ensuring that the population remains monomorphic 
between two events of invasion by a new trait.
\\

\textbf{Assumption "Invasion implies fixation" }
{\it For all} $x \in \mathcal{X}$ 
{\it and for almost every} $y\in \mathcal{X}~:$

{\it either} $\displaystyle \frac{b(y)-d(y)}{\eta(y)C(y-x)} < \displaystyle \frac{b(x)-d(x)}{\eta(x)C(0)},$
$\qquad $

{\it or } $\,\displaystyle \frac{b(y)-d(y)}{\eta(y)C(y-x)} > \displaystyle \frac{b(x)-d(x)}{\eta(x)C(0)}$ {\it and} $\displaystyle \frac{b(x)-d(x)}{\eta(x)C(x-y)} < \displaystyle \frac{b(y)-d(y)}{\eta(y)C(0)}.$
\\

\textbf{Remark : }
 $\hat{\mathrm{n}}_{x} := \textstyle \frac{b(x)-d(x)}{\eta(x)C(0)}$
 is the equilibrium of the dynamical system 
 that a large population size of individuals with trait $x$ approximate.

Moreover, in the case of logistic populations with a constant $C$, 
 this assumption is satisfied as soon as
 $x \rightarrow \hat{\mathrm{n}}_{x}$  is strictly monotonous.

\subsection{Flemming-Viot process and marker evolution}

The dynamics of the marker is first defined for a constant trait $x$ 
via a Flemming-Viot process as defined below.
This process generalizes to a potential infinity of markers
the measure-valued process describing two markers :
$F_t(dv) = X_t\, \delta_0(dv) + (1-X_t)\, \delta_1(dv)$
with $X$ 
the neutral Wright-Fisher diffusion 
defined in Section \ref{sec:aure}, cf. (\ref{WFdiff})
with $s = 0$ and remarks below.

In the sequel,
 we denote by $\mathcal{P}(\mathcal{U})$ 
 and $\mathcal{P}(\mathcal{X}\times \mathcal{U})$ 
 the probability measure spaces respectively
on $\mathcal{U}$ and on $\mathcal{X}\times \mathcal{U},$
while $\langle \nu |\ \phi\rangle$
denotes the integral of the measurable function $\phi$ against the measure $\nu$.
Also, $\mathcal{M}_F(\mathcal{X} \times \mathcal{U})$ 
is the set of finite measures on $\mathcal{X} \times \mathcal{U}$.

\begin{Def}
{\it Given } $x \in \mathcal{X}$ {\it and} $u \in \mathcal{U}$, {\it the Fleming-Viot process} $(F_{t}^{u}(x, t\ \geq\ 0)$
{\it indexed by} $x$, 
{\it started at time} $0$ {\it with initial condition} $\delta_{u}$ 
{\it and associated with the mutation operator} $A$ 
{\it is the} $\mathcal{P}(\mathcal{U})$-{\it valued process whose law is characterized as the unique solution of the following martingale problem. 
For any} $\phi\in \mathcal{D}(A)$ ,
\begin{align}
M_{t}^{x}(\displaystyle \phi)
=\langle F_{t}^{u}(x,\,.) |\, \phi\rangle-\phi(u)
-b(x)\, \textstyle \int_{0}^{t}\langle F_{s}^{u}(x,\,.) |\, A\phi\rangle\, ds
\label{FVeq}
\end{align}
{\it is a continuous square integrable martingale with quadratic variation process}
\begin{align*}
\langle &M^{x}(\phi)\rangle_{t}
=\displaystyle \frac{2b(x)}{\hat{n}_{x}}
\int_{0}^{t}(\langle F_{s}^{u}(x,\,.) |\, \phi^{2}\rangle
-\langle F_{s}^{u}(x,\,.) |\, \phi\rangle^{2})ds.
\end{align*}
\end{Def}

\textbf{Remark : }
The model presented in Section \ref{sec:aure}
provides simple illustrations for such kind of processes.
For the solution $X$ of (\ref{WFdiff}),
define :
$$
F_t(dv) := X_t\, \delta_0(dv) + (1-X_t)\, \delta_1(dv).
$$
Then, with $A \equiv 0$ 
(the only transitions for the traits are between 0 and 1),
$S$ a continuous function such that $S(0) = s$ 
(the parameter of selection at individual level) 
and $S(1) = 0$, 
for any $\varphi$ measurable :
\begin{align*}
&M_t(\phi) 
:= \langle F_{t} |\, \phi\rangle - \langle F_{0} |\, \phi\rangle
- (\langle F_{t} |\, S\times \phi\rangle - 
\langle F_{t} |\, S\rangle\times \langle F_{t} |\, \phi\rangle)
\\&\textstyle
= [\phi(0) - \phi(1)]\, (X_t - x - \int_0^t s\, X_r (1-X_r)\, dr)
= [\phi(0) - \phi(1)]\, \sigma\,  \int_0^t \sqrt{X_r (1-X_r)}\, dB_r,
\end{align*}
where we used Ito's formula, 
is clearly a square-integrable martingale with quadratic variation :
\begin{align*}
\langle M(\phi)\rangle_t 
=  [\phi(0) - \phi(1)]^2 \,  \sigma^2\, \int_0^t X_r (1-X_r) dr
= \sigma^2\, \int_0^t (\langle F_{r} |\, \phi^2\rangle - \langle F_{r} |\, \phi\rangle^2)\, dr.
\end{align*}
Note that $s$ is here a selective effect on the distribution
that is not present in (\ref{FVeq}) because the marker is neutral.
\\

Moreover, 
recalling the equation that described 
the state $\mu_t$ of the population of groups 
in Section \ref{sec:aure},
one can relate it to equation (\ref{FVeq})
with $b(x)\, A = \mathcal{L}_{W\!F}$
and the martingale $M(\phi)$ being identically zero
(i.e. with zero quadratic variations).
In \cite{LM15},
the authors in fact derive another description of the population of groups
in the limit of large population sizes (intra-groups and inter-groups).
This limit is also described as such Flemming-Viot process,
with a non-zero martingale 
because one does no longer neglect 
the non-selective
birth and death events of groups.
Again, 
its quadratic variations satisfies :
$d\langle M(\phi)\rangle_t 
\propto \int_{0}^{t}(\langle F_{u} |\, \phi^{2}\rangle
-\langle F_{u} |\, \phi\rangle^{2})\, du.$
Of course, 
there is still the additional term involving $r$ 
in the finite variation part (\ref{FVeq}).
I refer for instance to \cite{E00} and \cite{D96} 
for a detailed presentation of Flemming-Viot processes.
\\

\subsection{Convergence to the Substitution Flemming-Viot Process}
We can now state the main theorem 
that describes the slow-fast dynamics of adaptive traits
and neutral markers at the (trait) evolutionary time scale.

\begin{theo}
 We work under Definition \ref{DefT}
   and the assumption that "invasion implies fixation". 
Consider the initial condition
$\nu_{0}^{K}(${\it dy}, $dv) = n_{0}^{K}\delta_{(x_{0},u_{0})}(${\it dy}, $dv)$ 
{\it with} $\textstyle \lim_{K\rightarrow\infty}n_{0}^{K}
=\hat{\mathrm{n}}_{x_{0}}$ 
{\it and} $\textstyle \sup_{K\in \mathbb{N}^{*}}\mathbb{E}((n_{0}^{K})^{3}) 
<+\infty.$
{\it Then, the population process} $(\nu_{Kt}^{K},\ t \geq 0)$ {\it converges in law to the} $\mathcal{M}_{F}(\mathcal{X}\ \times \mathcal{U})$-{\it valued}
{\it process} $(V_{t}(dy, dv),\ t\geq 0)$.

To define this Markov process,
with initial condition $\hat{n}_{x_0}\, \delta_{x_0}\, \delta_{u_0}$,
we only need to describe it until the first jump of the trait,
which is given by an exponential law.
Namely, 
the trait jumps from $x_0$ to $x_0 + k$ with rate :
\begin{align*}
b(x_0)\, \displaystyle \hat{\mathrm{n}}_{x_0} \frac{[f(x_0+k;x_0)]_{+}} {b(x_0+k)}
m(x_0,\, k) dk.  
\end{align*}
Then, given that this first jump occurs at time $T$,
the law of the new marker is given by :
\begin{align*}
U \sim F_{T}^{u_0}(x_0, du),
\text{ so that }
V_T(dy,\, dv) = \hat{n}_{x_0+k}\, \delta_{x_0+k}\, \delta_{U}.
\end{align*}
Then, for any $t< T$, we have :
\begin{center}
$V_{t} (dy, dv) 
= \hat{\mathrm{n}}_{x_0}\delta_{x_0}(dy) F_{t}^{u_0}(x_0, dv)$,  
\end{center}
and the process is defined recursively like this before the next jumps.

The convergence holds in the sense of finite dimensional distributions 
on $\mathcal{M}_F(\mathcal{X} \times \mathcal{U})$.
In addition, 
the convergence also holds in the sense of occupation measures, 
i.e. the measure $\nu^K_{K T}(dy ;\, dv)\,dt$ on $\mathcal{X} \times \mathcal{U} \times [0; T]$ converges weakly to the measure $V_t(dy ;\, dv)\, dt$
for any $T > 0$.
\end{theo} 

This process is called by the authors the Substitution Fleming-Viot Process (SFVP).
It generalizes the Trait Substitution Process (TSP) introduced in \cite{M96}
and also obtained  as a limit of individual-based model in \cite{C06}.
The TSP is in fact the marginal of the SFVP on the trait space.
The jump rate of the TSP 
can be  easily interpreted.
The mutation rate 
of any resident 
is given by $b(x_0)\, m(x_0,\, k) dk$.
The number of such residents at equilibrium
is nearly $\hat{n}_{x_0}\, K$.
While in competition with the resident population,
the survival of 
the lineage of the mutant 
depends mainly on the period 
where the associated sub-population 
is too small to disrupt 
the resident population.
From classical results of Branching process,
it survives with probability
$[f(x_0+k;x_0)]_{+}/b(x_0+k)$.
Note that only beneficial mutations pass through,
so that they shall invade quickly 
after reaching a non-negligible proportion
in the population.
The product provides the rate of occurrence 
of such successful mutation in the whole population,
that gets divided by $K$ in the new time-scale.
\\

Given the recent and impressive progress
in sequencing and comparing 
genetic data between species,
one has partly access to the marker dynamics.
The selective dynamics is however much more difficult to follow,
since one would have
to evaluate 
mutation effects 
and the advantage they bring
in the past eco-systems.
It would thus be of high interest 
to be able to infer strong selective effects 
from the dynamics of the marker.

At each sweep, a very specific marker is selected.
This effect is referred to as a genetic hitchhiking.
It shall be chosen according to the law 
$F_{T}^{u_0}(x_0, du)$.
We may expect that numerous selective sweeps 
should increase the variability of the markers.
Yet, without assuming any effect on the marker dynamics 
from the selective component,
it is not very clear 
that doing these steps more frequently 
shall increase the variance in the selected marker
(after a comparable time span).

Nonetheless,
such hitchhiking events usually leave a mark,
at least for sexual reproduction species.
Because of frequent recombination events,
hitchhiking effect is mostly effective 
for markers closely linked to
the selected allele.
The diversity of variants 
become small
that one gets closer in the genome to the selected allele.
A well-known example is given by the selection for genes 
that favor the digestion of milk.
As one can imagine, 
the analysis is then
much more demanding 
than the convergence to the SFVP.

\subsection{The last time-scale of adaptive-dynamics}
\label{AD}

As the last time-scale presented 
in the current paper,
Adaptive Dynamics is probably the one that gives to natural selection
the most predictable effect.
Let us follow Baar, Bovier and Champagnat \cite{BBC17} :
the connection of Adaptive dynamics
 with the individual-based model
is demonstrated through a single step of convergence
(as long as no singularity is reached).
In this context,
the canonical equation of adaptive dynamics (CEAD)
states that the population of interest
can be considered monomorphic,
and its trait $x_t$ follows an ordinary differential equation
 of first order.
 Namely,
the speed of the trait involves
 the mutation rate,
 the population size at equilibrium,
some derivative of the fitness of invasion,
and the squared effect of mutations in the direction 
of invasion.
Remark that this last term 
is not the mean effect of mutations in this direction
because the more effective is a mutation, the more likely it is to fix.
In order to obtain such deterministic behavior,
we again need the assumptions for the TSS, 
i.e. rare mutations as compared to the ecological time-scale,
with negligible fluctuations around the size equilibrium
and invasion implying fixation
(cf. previous Subsection).
Moreover, the CEAD relies on the assumption 
that mutations have infinitesimal effects,
so that the trait evolves continuously
by the accumulation of
large number of  such mutations.
The connection of the TSS to the CEAD 
involves a coupled rescaling of time 
and of fitness effects,
which is rather natural.
Yet,
for the actual connection with the individual-based model,
it introduces the major difficulty 
that any favorable mutation step shall be quite insignificant
and yet shall replace effectively the dominant trait.
In the same idea,
the dominant population
shall prevent deleterious yet almost neutral mutations to invade
and filter favorable mutation with an invasion success
 still proportional to the mutation effect.
 
Nonetheless,
the authors of \cite{BBC17} actually manage 
to demonstrate a regime of convergence to the CEAD,
where these issues are rigorously controlled.
Their individual-based model is quite close to the previous one, 
except that there is no marker anymore
and that the possible mutation steps are assumed 
to be on some discrete and finite grid 
(whose mesh size goes to 0),
preventing large mutations.
I thus use the same notations as previously 
(rather than the one of \cite{BBC17})
and refer to Subsection \ref{sec:iBm}.
The main difference is also a scaling parameter
 for the mutation effect~:
 
{\it When a trait mutation occurs 
(in a population with trait $x$), 
 the new trait of the descendant is} $x+\sigma_K\, k \in \mathcal{X}$ 
 {\it with} $k$ {\it chosen according to the probability measure} 
 $\{m(x, k)\}_{k\in [\![-A, A]\!]}$.
 {\it The mutation rate is also allowed to depend on $x$,
 and is thus given by $q_K\, M(x).$}
 \\
 
 Besides the other assumptions 
 we have made in Subsection \ref{sec:iBm},
there are additional issues of regularity for 
the birth rate $b$, the death rate $d$, 
the mutation rate $M$,
the sensibility to competition $\eta$ 
and the competition kernel $C$,
for which assumptions are required. 
It is also assumed that $b(x) > d(x)$ 
and  that $C(x, x)$ is uniformly upper-bounded
for any $x\in \mathcal{X}$.
 
There is a last assumption
to ensure the absence of singularity,
based on the invasion fitness
of a mutant $y$ in a resident population $x$ :
$$
f(y, x) := b(y) - d(y) - \eta(y)\, C(y, x)\, \hat{n}_x.
$$
It indicates the mean growth rate of a mutant population with trait $y$ 
as long as it is still negligible as compared to the resident population with trait $x$.
As stated in the convergence to the TSS,
the invasion probability 
of a mutant population
initiated by a single individual with trait $y$
tends to $f(y, x)_+ / b(y)$
(for large population size).
Here, $f_+$ is the positive part of $f$,
meaning that deleterious mutations cannot invade.
\\
  
\textbf{  Assumption 3 :}
  For all $x\in \mathcal{X},  \qquad \partial_{1}f (x, x)\neq 0.$
  \\
  
  Assumption 3 implies that either $\forall x\in \mathcal{X}:\partial_{1}f (x, x)>0$ or $\forall x\in \mathcal{X}:\partial_{1}f(x, x)<$
  $0$. 
  Therefore, 
  coexistence of two traits is not possible. 
  Without loss of generality,
  we can assume that, 
  $\forall x\in \mathcal{X}, \partial_{1}f (x, x)>0$. 
  
\begin{theo}
\label{thm:AD}
 {\it Assume that Assumptions} 1 {\it and} 3 {\it hold and that there exists a}
  {\it small} $\alpha>0$ {\it such that}
  \begin{align*}
&    K^{-1/2+\alpha}\ll \sigma_{K} \ll 1
\\   \text{ and} \qquad
&   \displaystyle \exp(-K^{\alpha})\ll q_{K} \ll 
  \frac{\sigma_{K}^{1+\alpha}}{K\, \ln K}, \quad \text{ as } K\rightarrow\infty.
  \end{align*}
  {\it Fix} $x_0\in \mathcal{X}$ {\it and let} $(N_{0}^{K})_{K\geq 0}$ {\it be a sequence of} $\mathbb{N}$-{\it valued random variables such}  
  {\it that} $N_{0}^{K}/K$ {\it converges in law, as} $ K\rightarrow\infty$, {\it to the positive constant} $\hat{\mathrm{n}}(x_0)$ {\it and is}
  {\it bounded in} $\mathrm{L}^{p}$, {\it for some} $p>1.$
  
  {\it For each} $K\geq 0$, {\it let} $\nu_{t}^{K}$ {\it be the process generated by} $\mathcal{L}^{K}$ {\it with monomorphic}
  {\it initial state} $(N_0^K/K).\,\delta_{\{x_0\}}$. 
  {\it Then, for all} $T>0$, {\it the sequence of rescaled processes},
  $(\nu_{t/(K\, q_K\, \sigma^2_K)}^{K})_{0\leq t\leq T}$, 
  {\it converges in probability, 
  as} $ K\rightarrow\infty$, 
  {\it with respect to the Skorokhod topology on} 
  $\mathbb{D} ([0, T], \mathrm{M}(\mathcal{X}))$ 
  {\it to the measure-valued process} $\hat{\mathrm{n}}(x_{t})\, \delta_{x_{t}},$
  {\it where} $(x_{t})_{0\leq t\leq T}$ 
  {\it is given as a solution of the CEAD},
  \begin{center}
    $\displaystyle \frac{dx_{t}}{dt}
  = \sum_{k = -A}^{A} k\, [k\, M(x_{t})\, 
   \hat{\mathrm{n}}(x_{t})\,
   \partial_{1} f(x_{t}, x_{t})\, ]_{+}
   m(x_{t}, k),\qquad $
   \text{ {\it with initial condition} $x_0.$}
  \end{center}
\end{theo}

\textbf{  Remarks : }
  
  (i) The main result of the paper actually holds under weaker assumptions.
  More precisely, Assumption 3 can be replaced by the following~:
 
\noindent \textbf{ Assumption $3'$}. {\it The initial state $v_{0}^{K}$ has a.s. support $\{x_0\}$
  with $x_0\in \mathcal{X}$ satisfying $\partial_{1}f(x_0, x_0)\neq 0.$}
  \\
  
  The reason is that,
   since $x\mapsto\partial_{1}f (x, x)$ is continuous, 
    Assumption 3   is satisfied locally. 
Since moreover $x\mapsto\partial_{1}f (x, x)$ is Lipschitz-continuous, 
the CEAD   never reaches in finite time an evolutionary singularity 
  (i.e. a value $y\in \mathcal{X}$ such that $\partial_{1}f(y, y)=0$).
   In particular, 
   for a fixed $T>0$, 
   the CEAD only visits traits 
   in some interval $I$ of $\mathcal{X}$ 
  where $\partial_{1}f (x, x)\neq 0$. 
  By modifying the parameters of the model out of $I$ 
  in such a way that $\partial_{1}f (x, x)\neq 0$ 
  everywhere in $\mathcal{X}$, 
  we can apply Theorem \ref{thm:AD} 
  to this modified process $\tilde{\nu}$.
  Then, we deduce that $\tilde{\nu}_{t/(K u_{K}\sigma_{K}^{2})}$ 
  has support   included in $I$ 
  for $t\in[0, T]$ with high probability, 
  and hence coincides with $\nu_{t/(K u_{K}\sigma_{K}^{2})}$
  on this time interval.
  
  (ii) The condition $q_{K}\displaystyle \ll\frac{\sigma_{K}^{\mathrm{l}+\alpha}}{K\ln K}$ allows mutation events during an invasion phase of a mutant trait, 
  but ensures that there is no "successful mutation"
  event during this phase.
  
  (iii) The fluctuations of the resident population
   are of order $K^{-1/2}$,
   thus $ K^{-1/2+\alpha}<< \sigma_{K}$
   ensures that the sign of the initial growth rate
    is not influenced by the fluctuations of the population size. 
    If a mutant trait $y$ 
    appears in a monomorphic population with trait $x$, 
    then its initial growth rate is
  $b -d(y)-\eta(y)\,C(y, x)\langle\nu_{t}^{K}\,\vert\, 1\rangle
  =f(y, x) + o(\sigma_{K})
  =(y-x)\, \partial_{1}f(x, x)+o(\sigma_{K})$ since
  $y-x=O(\sigma_{K})$.
  
  (iv) $\exp(K^{\alpha})$ is the time 
  during which the resident population stays with high probability in
  an $O(\varepsilon\sigma_{K})$-neighborhood of an attractive domain. This is a moderate deviation
  result. Thus, the condition $\exp(-K^{\alpha})\ll u_{K}$ ensures that the resident population
  is still in this neighborhood when a mutant occurs.
  
  (vi) The time scale is $(K\, q_{K}\sigma_{K}^{2})^{-1}$ since the expected time for a mutation event
  is $(K \, q_{K})^{-1}$, the probability that a mutant invades is of order $\sigma_{K}$ and one needs
  $O(\sigma_{K}^{-1})$ mutant invasions to see an $O(1)$ change of the resident trait value. 
  \\

Still, 
such a strong filtering of mutations 
is probably the most questionable issue of realism
  concerning the modeling of evolution.
The fluctuations around the deterministic system 
   shall be extremely small 
   and slightly deleterious mutations
shall be well-separated for such conclusions to be satisfied.
   It seems unlikely that selective effects are so dominant
   even for rather large populations (a million of individuals or so).
Especially since there is usually a structuring 
of the population in term of non-heritable 
or loosely heritable characteristics,
that are not neutral for survival. 
Think for instance of the individual positions
or their level of infection by parasites.
As we have seen in Section 1,
it may induce much more variability 
compared to the case where all the individuals are identical.
  In this time-scale of infinitesimal mutations, 
  we may expect to see, 
  in addition to these selective effect 
  mainly driven by favorable mutations, 
  also some noise due to the fixation of almost neutral mutations.
  The trait of the population is still quite likely 
  to follow the direction given by invasion fitness,
  yet its displacement might be quite different 
  from the one given by the CEAD
  and not as regular.
  
\section{Conclusion}

As we have seen,
there is a large class of processes that 
can be rigorously obtained 
as limits of individual-based models
appropriately rescaled.
This is to be expected
since this representation fits the closest 
to simulations of populations,
with the minimal set of assumptions 
to include any interaction of interest.
Nonetheless,
the proof of convergence results are quite challenging
and impose to be very specific on the way
time-scales are well-separated.
By the coupled observation both of the proofs
and the simulations,
the main weaknesses of the models
usually appear much more salient.

Notably, 
we have evaluated the difficulty 
in estimating the birth rate
in Section \ref{sec:aline}
from the sole knowledge of the trait at birth 
(because there is a lot of fluctuations 
until the birth event).
In Sections \ref{sec:aure} and \ref{sec:henry},
the main issue appears to be that 
the predicted selection effects 
might be driven by too exceptional realizations
of the stochastic process 
describing the dynamics of a typical individual.
In simulations and actual life,
such transitions from a very stable equilibrium 
to another one quite separated
shall not happen exactly the way these models predict,
notably for the time at which they occur.
Although corrections can be made by some truncation,
or by adding another term to the equation 
governing the density,
the most realistic first step
would possibly be to consider the Flemming-Viot representation 
introduced in Section \ref{sec:tran}.
Yet, such measured-valued stochastic process
is quite more challenging to describe.
In Section \ref{sec:tran},
the assumption 
that mutations have only an infinitesimal effect
appears difficult to combine 
with the fact
that the population stays monomorphic
and that the mutations
are filtered depending on their effects.

This is the core of science 
to start with the most elementary models,
like the ODE defining the growth rate of a population,
and then to progressively incorporate 
more realistic features.
With the current probabilistic tools at our disposal,
it is clearly time 
to relate most of these models 
describing the dynamics of densities
to the individual-based measure-valued processes.
The main requirement is clearly 
that one averages over a large number of individuals,
but this can be obtained in very various ways 
depending on the interactions of interest.
In the case where a Central Limit Theorem holds,
it can be exploited 
to confirm the stability 
of the less noisy dynamics.
It may also provide
another dynamics,
a priori closer to IBM,
with different qualitative properties
as in Section 2.
Moreover,
the convergence results can be stated 
for very diverse time-scales,
from the rapid adaptation of cells 
to the propagation of parasites
and the evolution of species over millions of years.
They provide an elementary way 
to unify the models 
of micro-biology,
ecology and evolution.
Thus,
 they enable to justify more rigorously
 the separation of the related time-scales
or on the contrary motivate interesting couplings.
Notably, in Section \ref{sec:aure}  
the selective effects
 are closely linked
 to random demographic fluctuations;
 while in Section \ref{sec:tran},
 the evolution of a marker as a measure-valued process
 is coupled in a very specific way to the punctual 
 events of fixation for the traits under selection.
\\

\textbf{Acknowledgment}
\\
I wish to thank Etienne Pardoux,
Viet Chi Tran
and an anonymous reviewer
 for their contribution in improving 
 the readability of this article.
My thanks go also to Marc Hoffman 
and Beno\^it Henry 
for their answer.


\begin{thebibliography}{99}
\bibitem{BB18}
 Baar, M., Bovier, A.; 
 The polymorphic evolution sequence for populations with phenotypic plasticity,
    Electron. J. Probab.,
    Volume 23, no. 72, pp. 1-27 (2018)



\bibitem{BBC17}
Baar, M., Bovier, A., Champagnat, N.; 
From stochastic, individual-based models 
to the canonical equation of adaptive dynamics in one step. 
Ann. Appl. Probab. 27, no. 2, 1093--1170. (2017)

\bibitem{BDMT11} Bansaye, V., Delmas,J-F., Marsalle, L., Tran, V.;
Limit theorems for Markov processes indexed by continuous time Galton–Watson trees
    Ann. Appl. Probab.,
    V 21, Nbr 6, 2263-2314 (2011)

\bibitem{BM15}
Bansaye, V., M\'{e}l\'{e}ard, S.;
 Stochastic Models for Structured Populations, Scaling Limits and Long Time Behavior,
 Stochastics in Biological Systems, 1.4 (2015)

\bibitem{BEV17}
Barton N.H., Etheridge A.M., Véber A.;
The infinitesimal model: Definition, derivation, and implications.
Theor Popul Biol.;V 118, pp. 50-73 (2017) 

\bibitem{BCFMT}
Billiard, S., Collet,P., Ferri\`ere, R., M\'{e}l\'{e}ard,S., Tran, V.; 
The effect of competition and horizontal trait inheritance on invasion, fixation and polymorphism. 
Journal of Theoretical Biology, Elsevier, 411, pp.48-58 (2016)


\bibitem{bfmt18}
Billiard, S., Ferri\`ere, R., M\'{e}l\'{e}ard, S., Tran, V.; 
Stochastic dynamics of adaptive trait and neutral
marker driven by eco-evolutionary feedbacks
 J. Math. Biol.  71: 1211-1242 (2015)


\bibitem{C06} Champagnat, N.;
A microscopic interpretation for adaptative dynamics trait substitution sequence model. Stoch Process Appl 116:1127–1160 (2006)

\bibitem{CFM06}
Champagnat, N., Ferri\`ere, R., M\'{e}l\'{e}ard, S.;
 Unifying evolutionary dynamics: from individual stochastic processes to macroscopic models. Theoretical Population Biology, Elsevier, 69 (3),
pp.297-321 (2006)


\bibitem{Henry} 
Champagnat, N., Henry, B.; 
A probabilistic approach to Dirac concentration
in non-local models of adaptation with several resources,
Ann. Appl. Probab.,
V.29, N.4 (2019)


\bibitem{CM11}
Champagnat, N., M\'{e}l\'{e}ard, S.;  Polymorphic evolution sequence and evolutionary branching. Probab Theory Relat Fields 151:45–94 (2011)

\bibitem{D10} 
Dawson, D.A.; 
Introductory Lectures on Stochastic Population Systems,
Technical Report Series No. 451, Laboratory for Research in Stat. and Probab. (2010)

\bibitem{D96}  
Dawson, D. A.; 
Mesure-valued markov processes. 
 v. 1541 of Lectures Notes in Math., Ecole d'Et\'{e} de probabilit\'{e}s de Saint-Flour XXI, Springer, pp 1–260, New York (1993)

\bibitem{DJMP05}
 Diekmann, O.,  Jabin, P.E., Mischler, S., Perthame, B.;
  The dynamics of adaptation: An illuminating example and a Hamilton-Jacobi approach. Theor. Pop. Biol., 67:257–271 (2005)
  
\bibitem{Dur08} 
Durrett, R.; Probability Models for DNA Sequence Evolution, 2nd ed. Springer  (2008)
  
\bibitem{E00} 
Etheridge, A.. An introduction to superprocesses, University Lecture Series, v. 20. American
Mathematical Society, Providence (2000)




\bibitem{H99} 
Marc Hoffmann. Adaptive estimation in diffusion processes. Stochastic Process. Appl., 79(1):135-163  (1999)

\bibitem{Aline} 
Hoffmann, M., Marguet, A.;
Statistical estimation in a randomly structured branching population,
Stochastic Process. Appl., to appear https://doi.org/10.1016/j.spa.2019.02.015 (2019)

\bibitem{K13} 
Kutoyants, Y. A.; Statistical inference for ergodic diffusion processes. Springer Science and Business Media (2013)

\bibitem{LM15} 
Luo, S., Mattingly, J.; 
Scaling limits of a model for selection at two scales. 
Nonlinearity. 30. 10.1088/1361-6544/aa5499.  (2015)


\bibitem{MT09} 
M\'{e}l\'{e}ard, S., Tran, V.; 
 Trait Substitution Sequence process and Canonical Equation for age-structured populations.  J. Math. Biology, Springer (Germany), 58 (6), pp.881-921 (2009)

\bibitem{M96}
J.A.J. Metz and al. Adaptative
dynamics, a geometrical study of the consequences of nearly faithful reproduction. S.J. Van Strien and S.M. Verduyn Lunel (ed.), Stochastic and Spatial Structures of Dynamical Systems, 45:183–231 
(1996)

\bibitem{MBPS}
Mirrahimi S., Barles G., Perthame B., Souganidis P. E.: Singular hamilton-jacobi equation for the tail problem. SIAM J. Math. Anal. 44(6), 4297–4319 (2012)


\bibitem{PG09}  
Perthame, B., Gauduchon, M.; 
Survival thresholds and mortality rates in adaptive dynamics : conciliating deterministic and stochastic simulations. Math. Med. Biol. 27, 195–210 (2010)

\bibitem{S91}
Sznitman, A.S.; Topics in propagation of chaos. 
In Ecole d'Et\'{e} de Probabilit\'{e}s de Saint-Flour XIX—1989, Lecture Notes in Math.1464, pp 165–251. Springer, Berlin (1991).

\bibitem{V19}
Velleret, A.; Two level natural selection
under the light of quasi-stationary distributions,
on ArXiv : https://arxiv.org/abs/1903.10161
 (2019).
 
 \bibitem{WGGD06}
West, SA., Griffin, AS., Gardner A., Diggle SP.  Social evolution theory for microorganisms. Nat. Rev. Microbiol. 4:597–607 (2006)
\end{thebibliography}
\end{document}